\documentclass[11pt]{amsart}

\usepackage{amsmath,amssymb,amsfonts}

\usepackage{caption2}
\usepackage[]{hyperref}

\usepackage[pdftex]{color,graphicx}
\usepackage{multicol}
\usepackage{verbatim}

\numberwithin{equation}{section}

\setcaptionmargin{1cm}

\newcommand{\be}{\begin{equation}}
\newcommand{\ee}{\end{equation}}

\newcommand{\Li}{\text{Li}}

\usepackage{color}
\definecolor{rosso}{cmyk}{0,1,1,0.4}
\definecolor{rossos}{cmyk}{0,1,1,0.55}
\definecolor{rossoc}{cmyk}{0,1,1,0.2}

\begin{document}

\begin{center}
{
 \bf

Investigation about a statement  equivalent to Riemann Hypothesis (RH)
  }  
\vskip 0.2cm {\large

Giovanni  Lodone \footnote{
giolodone3@gmail.com}  %
.} 
\\[0.2cm]

{\bf

ABSTRACT    } 
\\[0.2cm]

\vskip 0.2cm {\large

\begin{comment}
We try to approach  a known equivalence to RH involving relative maxima and minima of $\xi(s)$ 
on critical line by a representation of the derivative of the phase of $\xi(s)$ with respect to imaginary coordinate that involves directly Euler product.  
In this attempt it is found an 
object conjectured to be the, almost everywhere,  converging  ``spectrum'' of prime numbers.  
Reasons and consequences of the conjecture are highlighted.
\end{comment}

%
First  idea   is  to compute a quantity like the angular momentum  with respect to $(0,0)$, of an unitary mass of coordinates  $(\Re [\xi(s)],\Im [\xi(s)])$  while $\Im[s]$ is the time, and, $\Re[s]= constant$. If we impose that the derivative along $\Re[s]$, at points $\Re[s]= 1/2$ is grater than zero, then, we find   exactly a known  RH equivalence statement about relative maxima and minima  of  $\xi(1/2+i \Im[s])$ along critical line.  After representing this fictitious  “angular momentum”  by Euler Product, and, using PNT as a tool, it can be proved that this positivity condition is granted   everywhere at least for      $\xi(1/2+i \Im[s]) \ne 0$. So, if the above equivalence is true, it is found that  off-critical line zeros must be excluded for $\zeta(s)$ function along all critical strip . Further analysis on Euler Product (Lemma 2) has evidenced others shorter ways to same objective. Besides the converging spectrum of prime numbers is highlighted as a  by-product.
 } 
\\[0.2cm]

\vskip 0.1cm

\end{center}

{

MSC-Class  :   11M06, 11M99

  } 
\vspace{0.1cm}

 {\it Keywords} :    Riemann hypothesis , Riemann spectrum of primes , Prime numbers , Riemann zeta function , Euler product

%%%%%%%%%%%%%%%%%%%%%%%

     \tableofcontents

                                           \section { Introduction}        
                                           
   \noindent This paper is  linked to a previous one \cite   {Giovanni Lodone Nov2024} on Dirichlet functions.                                                 
  
  \noindent We use the following known relations: Riemann $\zeta(s)$ function, 
 (see also  \cite  {Apostol:1976}, \cite  {Davenport1980}, 
\cite [p.~229]  {Edwards:1974cz}), is defined as the analytic continuation of the function:

\be \label {1p1}
\zeta(s)=\sum _{n=1}^\infty \frac{1}{n^s} = \prod_{\forall p} \frac{1}{1-\frac{1}{p^s}} \quad \mbox{ ($p$ prime)}  \quad \quad \Re(s)>1
\ee
The strip
  $0<\Re(s)<1 $ is called critical strip, the line $\Re(s)= 1/2$ is  the critical line,  and,  the infinite product is called Euler Product. The complex argument $s$ is expressed throughout as:
\be \label{1p2}
 s=\frac{ 1}{2}+\epsilon +it
 \ee
so that $\epsilon$ is the distance from critical line ,%
 and, $t$ is the imaginary coordinate.
 
\noindent A key function is $\xi(s)$ % \
  \cite[p.~16]{Edwards:1974cz}, that is real on the critical line and has the same zeros of $\zeta(s)$ in the critical strip:
  
\be \label {1p3}
\xi(s)=\Gamma \left( \frac{s}{2}  +1\right)(s-1)\pi^{-s/2}\zeta(s) = 
\left[ \frac{\Gamma \left( \frac{s}{2} +1 \right) }{ \pi^{s/2}} \right]  \ \times \   \left[  \zeta(s)(s-1) \   \right]
\ee

\noindent RH statement is: 

 \noindent{ \it `` the analytic continuation of $\zeta(s)$  has all  non trivial zeros,  $s_0$,
  on $\Re[s]=1/2$ ( i.e.  on critical-line) ``.}

 \noindent  The equivalence  stated in 
 \cite[p.~6] {Bombieri:2000cz} is reported here for easy reading:

\noindent  { \it ''The Riemann hypothesis is equivalent to the statement that all local maxima
of $  \xi(t)=\xi\left( \frac{1}{2} + i \ t  \right)$ are positive and all local minima are negative, . . .'' }. 

\noindent This article is structured as follows:
\begin{itemize}

\item In section 2  is shown that the above equivalence is linked to  $\frac{\partial \angle[\xi(s)]}{\partial t} $, and,a numerical example is given.  %

\item In section 3  a similar expression,  $\frac{\partial \Im \{\ln [\zeta(s)(s-1)]\} }{\partial t} $, is computed  using Euler product .  It is proved  in Lemma 2 that this representation  is valid at least for $\epsilon > 0$. %Besides reasons are listed to conjecture a more wide area of validity. % 

\noindent It is stressed  also that this result opens new shorter ways to prove RH.

\item In section 4, by application of  Lemma 2 , it  is proved, with similar arguments in   \cite  [p.~11] {Giovanni Lodone Nov2024} ,    that 

$$\forall t \ \ with \ \  \xi(t) \ne 0 \ \  we \  have \ \ \left(  \frac{d \xi(t)}{dt} \right)^2- \xi(t) \times \frac{d^2 \xi(t)}{dt^2} >0$$. 
 Where $\xi(t)= \xi(1/2+it)$.

\item In section 5 results are summarized and concluding remarks are given.

\item In appendixes some numerical example,  and some analytical development are given.

\end{itemize}

%%%%%%%%%%%%%%%%%%%%%%%%%%%%
%%%%%%%%%%%%%%%%%%%%%

\section {Analytic phase variation of $\xi(s)$ along t}

\noindent %     
The phase of  $\xi(s)$ function can be computed using  {\it Wolfram Mathematica}  by RiemannSiegelZ function, or, the result is the same  ( see \cite {Giovanni Lodone 2024}),  
 by  \cite {Giovanni Lodone 2021}.
 
 \noindent  Through a positive  scale factor $ F(t) =  \left(\pi/2 \right)^{0.25}t^{\frac{7}{4}} e^{- \frac{\pi}{4}  t}  $  %
 and  hyperbolic functions plus other terms here omitted.  The  expression is:
$$
\frac{-\xi(t,\epsilon)}{F(t)  e^{i\epsilon\frac{\pi}{4}}}  \approx  Z(t,\epsilon)=
2 \sum_{n=1}^N \frac{    \cosh \left[ \epsilon \  \ln\left( \sqrt{\frac {t}{2 \pi n^2}  }\right) \right]   }{\sqrt{n}}
\cos\left(   t   \ln\left(    \sqrt{\frac {t}{2  e \pi n^2} } \right) - \frac{\pi}{8}   \right) 
$$
\be + 2 i  \sum_{n=1}^N\frac{ \sinh\left[ \epsilon \  \ln\left( \sqrt{\frac {t}{2 \pi n^2}  }\right) \right]   }{\sqrt{n}}
\sin\left(   t  \  \ln\left(    \sqrt{\frac {t}{2  e \pi n^2} } \right) - \frac{\pi}{8}  \right )         +R_0(t,\epsilon) \label {ZSinhECosh1}
\ee
where $N= \left  \lfloor       \sqrt{\frac {t}{2   \pi } }   \right    \rfloor$ and: 
\be \label {R0} 
 R_0(t,\epsilon) = (-1)^{N-1} \left( \frac{2 \pi}{t}\right)^{1/4}  \left[     C_0(p)
 \right]
\ee
with $p =  \sqrt{\frac {t}{2   \pi } } -N$. In (\ref{R0}), the dependence from $\epsilon$ has been neglected  (\cite [p.~1]{Giovanni Lodone 2021}).
Besides (\cite [p.~16]{Giovanni Lodone 2021}):
\be \label {CZero}
C_0(0.5)%
\le C_0(p)= \frac{(\cos(2 \pi(p^2-p-\frac{1}{16}))}{\cos(2 \pi p)} \le \cos(\frac{\pi}{8}) 
=C_0(0)=C_0(1)
\ee

\begin{figure}[!htbp]
\begin{center}
\includegraphics[width=1.0\textwidth]{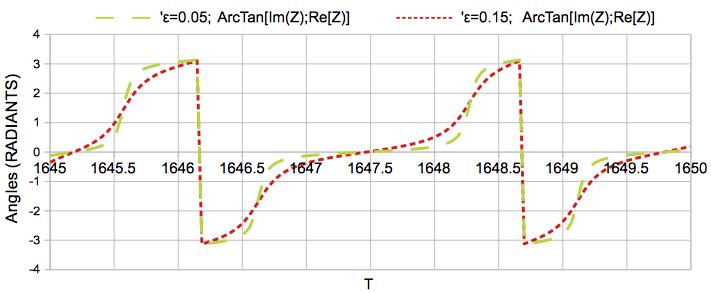} 
\caption{\small {\it    Phase of  $Z(t,\epsilon) = \angle [\xi(t,\epsilon)]$  ( \ref  {ZSinhECosh1} )
, in the interval  $ -\pi . .+\pi $ 
,    versus t
 , with $\epsilon$ as parameter. Zeros  are at  $t^*     \approx 1645.5737     \quad ,\quad t^*  \approx  1646.624 \quad  ,\quad t^*  \approx  1648.270 \quad  , \quad t^*  \approx 1649.118\quad $.
% ODS 17 02 21   sheet Eps = +0.05    i.e. 17_02_21_N5_8MT1645T1650
} }%
\label {AnglesVsTParamEps}  
\end{center}
\end{figure}

\begin{figure}[!htbp]
\begin{center}
\includegraphics[width=1.0\textwidth] {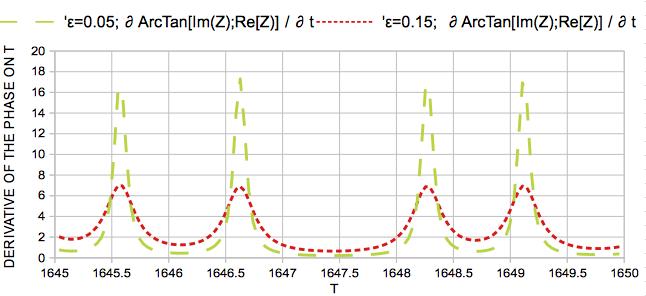} 
 \caption{\small {\it Derivative of $ \angle  [ \xi (t,\epsilon) ] $ % (
   with respect to $t$  by  phase of  $Z(t,\epsilon) $  ( \ref  {ZSinhECosh1} ),  %
  with $\epsilon$ as parameter.    %   % ODS 17 02 21 N5_8MT1645T1650.   sheet Eps005 N5 8M 
} }
\label {DerAnglesVsTParamEps}  
\end{center}
\end{figure}

    \noindent  We can compute  $\frac{\partial \angle[\xi(s)]}{\partial t} $ for $\epsilon \ne \epsilon^*_k$, where $\xi(1/2+\epsilon^*_k  + it^*_k)=0 \ \forall k$.  Of course for $\epsilon <  \epsilon^*$ the peaks of fig \ref {DerAnglesVsTParamEps}  computed with  $Z(t,\epsilon) $  are negative.  
We know that at least for  the first $10^{13}$ zeros $ \epsilon^*=0$, \cite  {Gourdon:2004cz}.
  \noindent While the phase variation of factor $\left[ \frac{\Gamma \left( \frac{s}{2} +1 \right) }{ \pi^{s/2}} \right] $ in (\ref {1p3}) is from  %
   \cite  [%p.~18-20
   appendix A] {Giovanni Lodone Nov2024}:

\be \label {DeFaseSuL1SuDeT}
\frac{\partial
\angle\left[ \Gamma\left(\frac{s}{2}+1\right)\pi^{- s/2} \right] }  {\partial t} = 
\left[ \frac { \partial [\Im_1+\Im_2+ \Im_3]} { \partial t } \right]=
   \frac {1}{2}\ln\left( \frac{t }{2 \pi}  \right)  + O(t^{-2}) \  ; \ \   for   \ \epsilon\rightarrow 0  \   then \ O(t^{-2}) \rightarrow  \frac{3  \epsilon}{4 t^2}
  \ee

 \noindent   For   $\ \epsilon\rightarrow 0  \ , \ O(t^{-2}) \rightarrow  \frac{3  \epsilon}{4 t^2} $ because with a few more additional developments ($ \left[ \frac { \partial^2 \Im_1} {\partial \epsilon \partial t } \right]_{\epsilon=0}=0$, and, neglecting $ \left[ \frac { \partial^2 \Im_3} {\partial \epsilon \partial t } \right]_{\epsilon=0}$) , we have:

 $$\left\{   \frac {\partial^2 \angle\left[ \left( \frac{\pi}{q}   \right)^{
 -\frac{s}{2}} \Gamma  \left( \frac{s+2}{2} \right)    \right]}{\partial \epsilon \partial t}  \right\}_{\epsilon=0} = \left[ \frac { \partial^2 [\Im_1+\Im_2+ \Im_3]} {\partial \epsilon \partial t } \right]_{\epsilon=0} 
 \approx
 \left[ \frac { \partial^2 \Im_2} {\partial \epsilon \partial t } \right]_{\epsilon=0}= \frac{3}{4t^2} $$

\noindent    The part of (\ref{1p3})  connected with peaks in fig. \ref {DerAnglesVsTParamEps}  at $t^*=\Im[\rho^*]$  is: %
 
     $$[ \zeta(s)(s-1) ] \ ,  \ i.e.  \   \  \frac{\partial \angle[ \zeta(s)(s-1) ]}{ \partial t}$$.

%%%%%%%%%%%%%%%%%%%%%%%%%%%%%%%%%%%%%%%%%%%%
\subsection {Angular Momentum of the $\xi$ function, and his  $\epsilon$ variation  } \
\label {ZeroOnCritLine}

\noindent As in
 \cite  [p.~3] {Giovanni Lodone Nov2024} we define angular momentum $\mathcal{L} [ \xi (s) ]$:

%$$ $$
 \noindent %
\be \label {AngMomDef}
DEFINITION \ 1 \ \ \ \   \ \ \ \ 
 \mathcal{L} [ \xi (s) ] :=
\det  \left(\begin{matrix} \Re\xi(s)  &  \Im\xi(s) \\ \frac{\partial}{\partial t}\Re\xi(s)  & \frac{\partial}{\partial t}\Im\xi(s)  \end{matrix}\right) 
=|\xi(s)|^2 \frac { \partial\angle[ \xi(s)]}{\partial t}
\ee

\noindent Since $\Im \xi=0$ for $\epsilon=0$, we have $\mathcal{L}[\xi(s)]=0$ on the critical line. From Lemma 1 in  \cite  [p.~3] {Giovanni Lodone Nov2024}  we have:  $\mathcal{L} [ C \times \xi (s) ]= |C|^2 \mathcal{L} [ \xi (s) ]$, where $C \in $ Complex field. So $Sign[\mathcal{L} [ \xi (s) ]]=Sign[\mathcal{L} [Z (t,\epsilon) ]$.

\subsection { LEMMA   1} \label {Lemma1} \

\noindent As in \cite  [p.~4, Lemma 3]  {Giovanni Lodone Nov2024}  we have:
\be \label {AngMomDerivSuEps_}
\ \ 
 \left[  \frac{\partial   \mathcal{L}  [ \xi (t,\epsilon) ]  )}{\partial \epsilon} \right]_{\epsilon=0}=
  det \left( \begin{matrix}
  \xi(t) &- \xi'(t )   \\ 
  \xi'(t) & -\xi''(t) \end{matrix}\right) %
   =[\xi'(t)]^2-\xi(t) \xi''(t) \ \ , \ where \ , \ \ \xi(t)) = \xi (1/2+it)
  \ee

    \noindent   PROOF: see \cite  [p.~4]  {Giovanni Lodone Nov2024}

\noindent It is apparent that equivalence  \cite[p.~6] {Bombieri:2000cz}, reported here p. 2, can be expressed with 

\be \label {StrongVersion}
 \left[  \frac{\partial   \mathcal{L}  [ \xi (t,\epsilon) ]  )}{\partial \epsilon} \right]_{\epsilon=0}=
[\xi'( 1/2+it)]^2-\xi(1/2 +it) \xi''(1/2+it) >0 \ \ \forall t
\ee

%\noindent We stress that, in this contest,  the statement : { \bf ``(\ref{StrongVersion}) is equivalent to RH''  is exactly 
\noindent So  this result  is linked to  equivalence  \cite[p.~6] {Bombieri:2000cz}, reported here p. 2.  The argument  is the same as in  \cite  [p.~4 Lemma3] {Giovanni Lodone Nov2024}.

  %

%

 %%%%%%%%%%%%%%%%%%%%%%%%%%%%%%
                                                                                                                                                                                                                                                                                                                                                                                                                                                                                                                                                                                                             \section {Phase variation of  $\xi(s)$ function using the Euler product} \label {DerivataAngLSuT}

                                                                                                                                                                                                                                                                                                                                                                                                                                                                                                                                                                                                       \noindent Following    \cite[p.~25]{Edwards:1974cz}))  we have:

\be \label {RootOfExplFormule_}
\ln(\zeta(s))=\ln(\xi(0))+ \sum_\rho \ln\left( 1-\frac{s}{\rho}\right)- \ln \left( \Gamma\left( \frac{s}{2} +1\right)  \right) + \frac{s}{2}\ln(\pi) - \ln(s-1) \ \ ; \  \{\Re(s)>1\}
\ee
                                                                                                                                                                                                                                                                                                                                                                                                                                                                                                                                                                                                               
                                                                                                                                                                                                                                                                                                                                                                                                                                                                                                                                                                                                               \noindent But notice that for at least $ \{\Re(s)>0\ \  , \forall t \ \ne t^*_k \} $  :

\be \label {RootOfExplFormule_2}
\ln(\zeta(s)(s-1))=\ln(\xi(0))+ \sum_\rho \ln\left( 1-\frac{s}{\rho}\right)- \ln \left( \Gamma\left( \frac{s}{2} +1\right)  \right) + \frac{s}{2}\ln(\pi)  \ \ ; \  \{\Re(s)>0\ \  , \forall t \ \ne t^*_k \} \ at \ least
\ee
                                                                                                                                                                                                                                                                                                                                                                                                                                                                                                                                                                                                            \noindent %  
                                                                                                                                                                                                                                                                                                                                                                                                                                                                                                                                                                                                            \noindent By (\ref  {RootOfExplFormule_}) Riemann found the distribution of  $J(x)$, %(
                                                                                                                                                                                                                                                                                                                                                                                                                                                                                                                                                                                                            \be \label{GeiX}
J(x)= \frac{1}{2} \left[ \sum_{ p^n<x} \frac{1}{n} + \sum_{ p^n\le x} \frac{1}{n}  \right]
\ee
, i.e. the ``number of prime powers till $x$'',  from  (\ref {RootOfExplFormule_}) for $x>1$ (  \cite[p.~34]{Edwards:1974cz}) :

\be \label  {PrimesPowerCount}
J(x)=Li(x)-\sum_\rho Li(x^\rho) -\ln(2)+\int_x^\infty \frac{dy}{y(y^2-1)\ln(y)}=Li(x)-\sum_\rho Li(x^\rho) -\ln(2) %
 \ \  %
 % \ee
 \ee
 Or  ( \cite[p.~36] {Edwards:1974cz}): 
  
  \be \label{fill1}
  dJ(x) = \left[ \frac{1}{\ln(x)}  -2 \sum_{\Im \rho>0} \frac{\ x^{\Re \rho-1}}{\ln(x)} \cos[\Im \rho \ln(x) ]
  -   \frac{1}{x(x^2-1)\ln(x)}
   \right] dx
    \ \ \ \ ;  \ \ \ \ x>1
  \ee
                                                                                                                                                                                                                                                                                                                                                                                                                                                                                                                                                                                                              
                                                                                                                                                                                                                                                                                                                                                                                                                                                                                                                                                                                                              \noindent The term $\int_x^\infty \frac{dy}{y(y^2-1)\ln(y)}$, wich stems from trivial zero, is a term that can  be neglected, and  $ \frac{1}{\ln(x)} $ is the  density of the local mean value of prime powers  . 
So :
                                                                                                                                                                                                                                                                                                                                                                                                                                                                                                                                                                                                              
   \be \label {DistribXGrande0_}
J(x) - Li(x)  = - \sum_\rho Li(x^\rho) - \ln(2) \quad \rightarrow \quad d[  J(x) - Li(x) ] =
 -2 \sum_\rho \frac{dx \ x^{\Re \rho-1}}{\ln(x)} \cos[\Im \rho \ln(x) ]; \quad \rho=\Re \rho\pm i \Im \rho
\ee

\noindent $J(x)$ is related to $\pi(x)$, the number of prime first powers till $x$,   by \cite[p.~33]{Edwards:1974cz}):  %
 \be \label {fill2}
 J(x) = \sum_{n=1}^{ \left \lfloor \frac{\ln(x)}{\ln(2)} \right \rfloor    } \frac{\pi(x^{1/n})}{n}
  \ \  \ , \  and \  the \  inverse , \  \ \  \pi(x)=\sum_{n=1}^{ \left \lfloor \frac{\ln(x)}{\ln(2)} \right \rfloor    } \frac{\mu(n)}{n} \ J(x^{1/n}) 
\ee
                                                                                                                                                                                                                                                                                                                                                                                                                                                                                                                                                                                                              \noindent where $\mu(n)$ is Moebious function \cite  {Apostol:1976}.
                                                                                                                                                                                                                                                                                                                                                                                                                                                                                                                                                                                                              
\noindent
 Riemann used a  transform relation between  $J(x)$ and $\zeta(s)$ \cite[p.~22]{Edwards:1974cz}) (  integral  is in Stieltjes sense \cite {Apostol:1974}  
                                                                                                                                                                                                                                                                                                                                                                                                                                                                                                                                                                                                              )
                                                                                                                                                                                                                                                                                                                                                                                                                                                                                                                                                                                                           giving  (\cite[p.~22]{Edwards:1974cz})):
 
\be \label {Trasf1_}
\ln(\zeta(s))=
\sum_p \sum_n \frac{1}{np^{n s}} =\int_0^\infty x^{-s} d[J(x)] \ \ ; \ \ \Re[s] >1
\ee
 
\noindent Whose inverse trasform, applied to all pieces of (\ref{RootOfExplFormule_}), (see tab. \ref {tabRiemannTrasf})
can be written as  ( \cite[p.~25]{Edwards:1974cz})) :

                                                                                                                                                                                                                                                                                                                                                                                                                                                                                                                                                                                                               \be \label {RiemIntByPart} J(x)= - \frac{1}{\ln(x) 2 \pi i}\int_{a-i\infty}^{a+i\infty}\frac{d}{ds}\left[  \frac{\ln(\zeta(s)}{s} \right] x^s ds \quad ; \quad a>1
                                                                                                                                                                                                                                                                                                                                                                                                                                                                                                                                                                                                               \ee

 \noindent % 
  We are interested  %
in  $\angle[\zeta(s)(s-1)]=\Im[\ln( \zeta(s)(s-1)]$  of ( \ref {RootOfExplFormule_2}), and,  in using Euler Product at $\Re(s)<1$.

  \noindent Notice that  condition $Re[s] >1$ in (\ref{Trasf1_}) is due  to left term: $\ln(\zeta(s))$, and to the purpose of using trasforms (see tab. \ref {tabRiemannTrasf} ).

 \noindent Notice also that $\angle[z]= \Im[\ln(z)]$, but caution is needed if $ z$ is real negative,
where $\arctan\left( \Im[z] , \Re[[z]\right)$ %
is defined, and, can be continuous   ,while complex logarithm experiences a $2 \pi$  phase discontinuity, so makes sense that $\frac{\partial \Im[\ln(\zeta(s)(s-1))]}{\partial t} $ diverges to $+\infty$  across $ \  \ t=t^* , \epsilon \le \epsilon^* \ \ $, with $\zeta(1/2+\epsilon^*+it^*) = 0$.

\noindent Let us focus on  left part of ( \ref {RootOfExplFormule_2}). It  converges for  $ \{\Re(s) >0 \} \ \  \smallsetminus \{ Singular \ p. \} $.  

\noindent For $ Re(s) >0$  singular points are: the $\zeta(s)$ non trivial  zeros, $s=\rho^*_k=\frac{1}{2} +\epsilon^*_k + i t^*_k$,   and the segments from them  normal to  axis $ Re(s) = 0$ where  there is  a  discontinuity,  along $t$, i.e. a $2 \pi$  step,   of $\angle[\zeta(s)(s-1)] $ computed as $\Im[ \ln(\zeta(s)(s-1))] $.

\noindent Because the two factors of  (\ref {1p3} ), %
are defined $\forall s$,  while  trasforms only for $\Re(s)>1$,  we can write, using 3rd row of  tab. \ref {tabRiemannTrasf}  , at least %(  
:
 \be \label {Trasf4_}
\ln[\zeta(s)(s-1)]=
\int_0^\infty x^{-s} d[J(x) -Li(x)] \ \ ; \  \ \ \ ; \ %
\ee  
   
 \noindent  Left term of (\ref  {Trasf4_}) is defined  surely for $\{\Re(s)>0\}  \ \  \smallsetminus \{ Singular \ p. \} $,   right term, let us say,  for  $\{\Re(s)> \alpha\}$, where $\alpha$ is  discussed below.

\begin{table}[bht]
\begin{center}
\begin{tabular}{c|c|c|c|c|c|}
 a     & $\phi(x) :  \ln(\phi(x))= \int_0^\infty x^{-s} d[\Phi(x)]$   &  $ \Phi(x) :  \Phi(x) = - \frac {\int_{a-i\infty}^{a+i\infty} \frac{d}{ds}\left[\frac{\ln(\phi(s)] }{s} \right]   x^s ds }{\ln(x) 2 \pi i}$ & reference  \cite   {Edwards:1974cz}     \\
\hline
$a >1$ & $\phi(s) = \zeta(s)$ &$\Phi(x)= J(x)$ & pp.22-25\\
\hline
$a >1$ &$\phi(s) =\frac{1 }{(s-1)}$    & $\Phi(x)= Li(x)$ & pp.26-29 \\
\hline
$a >1$ &$\phi(s) =\Pi_\rho\left(1-\frac{s }{\rho} \right)$    & $\Phi(x)=- \sum_{\Im(\rho)>0}[ Li(x^\rho)+Li(x^{1-\rho})  ]$  & pp.29-30  \\\hline
$a >1$ &$\phi(s) =\pi^{s/2}$    & $\Phi(x)=0$  & p.31  \\
\hline
$a >1$ &$\phi(s) =\xi(0)$    & $\Phi(x)= \ln(\xi(0))$& p.31  \\
\hline
$a >1$ &$\phi(s) =\Gamma \left( \frac{s}{2}+1\right)$    & $\Phi(x)=\int_x^\infty\frac{dy}{y(y^2-1)\ln(y)} $ & p.32  \\
\hline
\end{tabular}
\end{center}
\caption{\small Trasform pairs used by Riemann to find (\ref {PrimesPowerCount}).  %
}
\label {tabRiemannTrasf}
\end{table}

    \noindent  We want to extract  phase  variation of $\zeta(s) (s-1)$ along $t$ from (\ref {Trasf4_}).

 \noindent In order to trade resolution along $t$ with computational burden  it is convenient to address to mean phase variation in an interval $(t_1,t_2)$. 
   
\noindent %
We consider $t$   increment  in order to compute        incremental ratio with respect to $t$, %
 and, define the value $p^*$ so that:
\be \label {Estremi}
t_2=t+ \frac{\pi}{\ln(p^*)} \quad \quad  \mbox{ ; }  \quad \quad t_1=t - \frac{\pi}{\ln(p^*)}
 \mbox{ ; } \ \ \  \Delta t= \frac{2\pi}{\ln(p^*)}
\ee

    \noindent The phase  of $\zeta(s)$, through Euler product, is  given by:

\be \label {EulerPrPhaseA}
\angle [ \zeta_{EP}(s)] =\Im[\ln(\zeta_{EP}(s))] = - \sum_{p  }\angle \left [  1-\frac{1}{p^s} \right] %
=
   \sum_{p  } \arctan \left( \frac{-
 \sin(\ln(p) t      )}{ p^{1/2+\epsilon}-  
\cos( \ln(p) t   )  }\right) \ \, \ \epsilon >1/2
\ee

      and, using an {\it ad hoc} notation, mean phase variation through Euler Product :

 \be \label {EulPhaseDerSmoothing2}
\left \{ \frac {\Delta \angle [ \zeta_{EP}(s)] }{\Delta t} \right\}_{p^*,p_{max}}
= -  \left. \frac {\ln(p^*)} {2 \pi}   \sum_{p \rightarrow p_{max}   } \arctan \left( \frac{
 \sin(\ln(p) t      )}{ p^{1/2+\epsilon}-  
\cos( \ln(p) t   )  }\right)\right|_{t_1}^{t_2} \ \, \ \epsilon >1/2
  \ee

  \noindent Notice, like in  \cite  [p.~6]  {Giovanni Lodone Nov2024}  { \bf we do not interchange differentation with infinite sum}, and while $p^*$ is big but fixed,  instead, $p_{max} \rightarrow \infty$. We refer at this procedure, like in   \cite  [p.~6]  {Giovanni Lodone Nov2024},  with the shortcut :  $''p_{max}     \ge p^* \rightarrow \infty''$ that must not be taken literally (it means: for each $p^*$ big, $p_{max} $ grows without bounds). 
 \noindent  %
 The Imaginary part of $- \int_0^\infty x^{-s} d[Li(x)]$, %
 is  :

  \be \label {Trasf1ImPart}
  \Im   \left [- \int_0^\infty x^{-s} d[Li(x)] \right]  =
 \int_0^\infty \frac{\sin( \ln(x) t)}{x^{1/2+\epsilon}} d[Li(x)] \ \ \ , \ \epsilon >1/2
\ee

 \noindent  We want to find the mean value of phase variation in interval $( t_1,t_2)$, like in (\ref {Estremi}),  as above, so, we set $\Delta t = t_2-t_1$ and define accordingly  the value $p^*$. % 

 \noindent As $\sin(\ln(p)t_2)-\sin (\ln(p)t_1) =2\cos(\ln(p)t) \sin\left( \frac{\pi \ln(p)}{\ln(p^*)}\right)$, using above {\it ad hoc} notation, the mean phase variation  in  $t_1, t_2$ is:

  \be \label {Trasf1ImPartDerivSuT}
\left \{  \frac{ \Delta \left \{ \Im [  \int_{p=2}^{p=p_{max}} x^{-s} d[Li(x)] \right\}  }{\Delta t} \right \}_{p^*,p_{max}} = -  %
\frac{\ln(p^*)}{ \pi}
\int_{p=2}^{p=p_{max}} \frac{\cos( \ln(x) t)}{x^{1/2+\epsilon}}
\sin\left(  
\frac{\pi \ln(x)}{\ln(p^*)}
\right) d[Li(x)] \ \\ \ , \ \epsilon >1/2 %
\ee

  \noindent   Using incremental ratio for (\ref{DistribXGrande0_}) with  (\ref {Estremi}) as above, we have:

  $$  \left\{  \frac{
   \Delta \left \{ \Im [ 
  \int_{p=2}^{p=p_{max}} x^{-s} d[J(x) -Li(x)] 
   ] \right\} 
    }{\Delta t}  \right \}_{p^*,p_{max}}=
     \left\{ 
\frac{ \Delta \left \{ \Im [ - \int_{p=2}^{p=p_{max}} x^{-s} d[\sum_\rho Li(x^\rho)] \right\}  }{\Delta t}
\right \}_{p^*,p_{max}}=$$

       $$
   \frac{\ln(p^*)} { \pi}\int_2^{p_{max}} 
     \frac{\cos(\ln(x) \ t)}{x^{1/2+ \epsilon}} \sin\left(  \frac{\pi \ln(x)}{\ln(p^*)} \right)  
   d \left [ \sum_{\rho_k} Li(x^{\rho_k}) \right]   %
   \quad  =
$$

 $$
  \frac{\ln(p^*)} { \pi}\int_2^{p_{max}} 
     \frac{\cos(\ln(x) \ t)}{x^{1/2+ \epsilon}} \sin\left(  \frac{\pi \ln(x)}{\ln(p^*)} \right)  
     \left[          
     2 \sum_{\rho_k } \frac{dx \ x^{\Re \rho_k-1}}{\ln(x)} \cos[\Im \rho_k \ln(x) ]      \right]
     =
 $$

\be \label {2Integrals_}
\frac {\ln(p^*)} { \pi}
  \int_{2} ^ {p_{max}}  \sum_{\rho^* }  \frac{dx \ x^{\epsilon^*-\epsilon}}{x \ln(x)} \cos[\Im \rho \ln(x) ]     %
 \cos[ t \ln(x) ] \sin\left( \frac{\pi \ln(x)}{\ln(p^*)}  \right)  %
\ee
 
  \noindent  Where  : $\quad \rho^*_k=\Re \rho^*_k\pm i \Im \rho^*_k$ are the non trivial zeros. %, 
  If $\epsilon^*_{MAX}= Max [\epsilon^*_k] , \forall k$, then
  %
  %
% % 
 \noindent  { \bf  for $\epsilon >\epsilon^*_{MAX}$, (\ref {2Integrals_}) converges absolutely}   $\forall p_{max}  \ge p^* \rightarrow \infty$.

   \noindent %
   This is apparent posing  $e^y=x$, so for (\ref  {2Integrals_}) :

 \be \label{fill5}
   \left\{  \frac{
   \Delta \left \{ \Im [ 
  \int_{p=2}^{p=p_{max}} x^{-s} d[J(x) -Li(x)] 
   ] \right\} 
    }{\Delta t}  \right \}_{p^*,p_{max}}=
     \left\{ 
\frac{ \Delta \left \{ \Im [ - \int_{p=2}^{p=p_{max}} x^{-s} d[\sum_\rho Li(x^\rho)] \right\}  }{\Delta t}
\right \}_{p^*,p_{max}}=
\ee

 \be \label {fill4}
  \frac{\ln(p^*)} { \pi}\int_{\ln(2)}^{\ln(p_{max})} 
     \frac{\cos(y \ t)}{e^{(1/2+ \epsilon)y}} \sin\left(  \frac{\pi y}{\ln(p^*)} \right)  
     \left[          
     2 \sum_{\rho_k } \frac{dy \ e^{( \epsilon^*_k+1/2)y}}{y} \cos[\Im \rho_k  y]      \right]
     =
 \ee

\be \label {ZeroOffEAllaY_}
 \frac{\ln(p^*)}{ \pi}  \int_{\ln(2)} ^  {\ln(p_{max})}
\sum_{\rho^*} \frac{dy \ e^{y(\epsilon^*-\epsilon)}}{ y} \cos[\Im \rho y ]     %
 \cos[y t] \sin\left( \frac{\pi y}{\ln(p^*)}  \right); \quad \rho=\Re \rho\pm i \Im \rho \ ; \ \epsilon^* < 1/2
\ee
 \noindent  The square bracket in intermediate lines of  (\ref{2Integrals_}) and (\ref{ZeroOffEAllaY_} is $d(J(x)-Li(x))$. The factor  $ e^{(\epsilon^*_k-\epsilon)y}$ in  (\ref{ZeroOffEAllaY_} forces  the integrals to converge for $\epsilon> \epsilon^*_k$. So $\alpha$ of definition halfplane for  (\ref {Trasf4_}) right term  is at least 
 
 \be \label{fill3}
 0.5+ \epsilon^*_{MAX} = \Re(s) =\alpha  <1
 \ee

 \noindent because no  roots $\rho^*$ on $\Re(s)=1$ are possible. See \cite[p.~69] {Edwards:1974cz}).

\noindent %{
From previous Lemmas,  and considering third row of tab. \ref {tabRiemannTrasf}  , we can  write %
right part of (\ref  {Trasf4_}),  at least for $\Re(s) >1$, as:

      $$\left\{  \frac{ \Delta \Im [\ln \{ \zeta(s)(s-1) \}]}{\Delta t} \right\}_{p^*, p_{max}} = %
        \left \{ \frac {\Delta \angle [ \zeta_{EP}(s)] }{\Delta t} \right\}_{p^*,p_{max}}  - \left \{  \frac{ \Delta \left \{ \Im [  \int_{p=2}^{p=p_{max}} x^{-s} d[Li(x)] \right\}  }{\Delta t} \right \}_{p^*,p_{max}} =$$

       \be \label {ProVarPhaseCsi}
 -  \frac {\ln(p^*)} {2 \pi}   \left\{ 
 \left[ \sum_{p =2}^{p_{max}   } \arctan \left( \frac{
 \sin(\ln(p) t      )}{ p^{1/2+\epsilon}-  
\cos( \ln(p) t   )  }\right)\right]_{t-\frac{\pi}{\ln(p^*)}}^{t   +\frac{\pi}{\ln(p^*)}   } - 
 \int_2^  { p_{max} }  
  \frac{2\cos(\ln(x) \ t)}{x^{1/2+ \epsilon}} \sin\left(  \frac{\pi \ln(x)}{\ln(p^*)} \right) 
     d [ Li(x) ] \right\} =
\ee  

\be \label {FromJminusLi}
 = - \frac{\ln(p^*)} { \pi}\int_{2}^{p_{max} \rightarrow \infty} 
     \frac{\cos(\ln(x) \ t)}{x^{1/2+ \epsilon}} \sin\left(  \frac{\pi \ln(x)}{\ln(p^*)} \right)  
   d \left [J(x)-Li(x)) \right]  
\ee

  % {ProVarPhaseCsi}   {FromJminusLi}

   \noindent  Below we %speculate 
   see
   that  right part of (\ref  {Trasf4_}), i.e. like (\ref{FromJminusLi}) or (\ref{ProVarPhaseCsi}),  have actually converging  features  well beyond   $\Re[s]>1$.

%%%%%%%%%%%%%%%%%%%%%%%%%%%%%%%%%%%%%%%
%%%%%%%%%%%%%%%%%%%%%%%%%%%%%%%%%%%%%%%%%%% 
  \subsection { LEMMA 2 Euler product in critical strip, and, extended convergence of  (\ref {ProVarPhaseCsi}) and  (\ref {FromJminusLi}) 
  %, and, shorter ways toward RH 
  } \label {CONJ} \

%   \begin{comment}

  \noindent  
Let us consider now the Euler product  indipendently from $(\Re(s)>1)$. This means that the relation with $\zeta(s)$ in  (\ref  {1p1}  ), wich is based on absolute convergence of Euler product for   $\Re(s)>1$,  is not granted. We have:

\be \label {LogEP}
\ln\left[  \prod_{\forall p} \frac{1}{1-\frac{1}{p^s}} \right]  = - \sum_{\forall p} \ln\left[ 1-\frac{\cos(\ln(p) t )-  i \sin(\ln(p) t  )  }{p^{1/2+\epsilon}} \right]
\ee
\noindent Using the taylor expansion $\ln(1+x)=x-\frac{x^2}{2}+\frac{x^3}{3} - . . .$ for each  logarithm in (\ref{LogEP}) we have:

\be \label {LogEP2}
\ln\left[  \prod_{\forall p} \frac{1}{1-\frac{1   }{p^s}} \right]  = -
\sum_{\forall p} 
\frac{-\cos(\ln(p) t) + i \sin(\ln(p) t)    }{p^{1/2+\epsilon}}- 
\frac{\left(  \frac{ -\cos(\ln(p) t) +  i \sin(\ln(p) t)    }{p^{1/2+\epsilon}}  \right)^2} {2} +
\frac{\left(  \frac{ -\cos(\ln(p) t) +  i \sin(\ln(p) t)    }{p^{1/2+\epsilon}}  \right)^3} {3} - . . . .
\ee

%    If we are interested only on  convergence  higher orders are not important.
%

\noindent If we are interested  only in convergence  for $\epsilon>0$, we  could neglect all terms but the first in (\ref{LogEP2}):
\be \label {parteImportante} 
\sum_{\forall p} 
\frac{\cos(\ln(p) t) - i \sin(\ln(p) t)    }{p^{1/2+\epsilon}} \ee
 and, realize that the imaginary part after treatment with mean variation estimation as by (\ref{Estremi})  gives exactly the approximate  expression   %(\ref {FdaLgaritmicInt})  for imaginary part, i.e. 
  for the  phase variation along $t$ (\ref{Approx}).
  
   \noindent  For real part of (\ref{LogEP2}), i.e. the logarithm of $|E.P.|$, a similar expression holds.

\noindent  Because for $\epsilon>0$,  absolute value of (\ref{LogEP2}), except first sum (i.e.  erasing (\ref{parteImportante}) ), is upperbounded by 
$$ \ \  \   \sum_{n>1 }^\infty \ \sum_{\forall p}  \frac{1}{n p^{(1/2+ \epsilon)n}}  <  \sum_{n>1 }^\infty \ \int_2^\infty  \frac{dy}{ y^{(1/2+ \epsilon)n }}
=
\sum_{n>1 }^\infty \frac{1}{[(1/2+ \epsilon)n] \times  2^{(1/2+ \epsilon)n-1}   }<
$$
$$
%  \sum_{n>1 }^\infty \frac{1}{[(1/2+ \epsilon)n] \times  2^{(1/2+ \epsilon)n-1}   } < 
\sum_{n>1 }^\infty \frac{1}{ \epsilon \times 2^{(1/2+ \epsilon)n-1}   } 
=
\sum_{n>1 }^\infty \frac{2/\epsilon}{  2^{(1/2+ \epsilon)n}   } 
<
\sum_{n>1 }^\infty \frac{2/\epsilon}{  2^{( \epsilon)n}   }=
\frac {2/(\epsilon 2^{ 2 \epsilon} )} {1-2^{-\epsilon}}
<
\infty  \ \ for \ \epsilon >0$$

\subsection {Lemma  2A: convegence of  $\ln[ EP(s)] - \int_0^\infty x^{-s} d[ Li(x)]$} \label {conv2_12}

%%%%%%%%%%%%%%%%%%%%%%%%%%%%%%
 
 { \bf Expression  $\ln[ EP(s)] - \int_0^\infty x^{-s} d[ Li(x)]$ %(\ref {Diff})
 , irrespective of the meanning that could be assigned  to it, 
 % $\left\{  \frac{ \Delta \Im [\ln \{ \zeta(s)(s-1) \}]}{\Delta t} \right\}_{p^*, p_{max}} $ for $\epsilon>1/2$, 
 is convergent for  $\epsilon >0$}.
% }

 We will prove
  for example  this for $\Re[ \ln(EP(s))]$ .  For divergence study ( neglecting all terms of  (\ref{LogEP2}) but (\ref {parteImportante}) ): 
                                                                                                                                                                                                                                                                                             \be \label {Third}
                                                                                                                                                                                                                                                                                             \Re \left[
\ln[ EP(s)] - \int_0^\infty x^{-s} d[ Li(x)] \right]
\approx 
 \Re\left[ 
\left\{
 \sum_{\forall p} 
\frac{\cos(\ln(p) t)-  i \sin(\ln(p) t)    }{p^{1/2+\epsilon}}
\right\} 
-\int_0^\infty x^{-s} d[ Li(x)] 
\right] <\infty \  ;  \ \forall t \ , \  \forall \epsilon >0
 \ee
  \noindent PROOF

  \noindent We take only first term in (\ref{LogEP2}),  i.e. (\ref{parteImportante} ).  In each interval containing only the prime p  the discrete sum  can be written as a defined integral, and, the infinite integral in $y$ can be split in an infinite sum of defined integrals.

 \noindent Supposing the prime $p$ (fixed in each interval, while $t$ is fixed everywhere)  is surrounded by neighbours primes at $p- 2\alpha(p) \ln(p) < p <p + 2\beta(p) \ln(p) $ 
 \footnote{The $p$ dependence of $\alpha(p)$, and $\beta(p)$ is implied in the following.  If the  distance among primes was exactly the mean distance, then $\alpha=\beta=1/2$. This is not the case so $\alpha$  and $\beta$ manage the practical cases. Notice a potential application to $h-$class primes $mod \ q$ ( see \cite  {Giovanni Lodone Nov2024}). The interval with mean $\alpha,\beta = 1/2$ would become: 
 $$ p_{h-class} -\phi(q) \alpha \ln( p_{h-class} ) < y <p_{h-class} +\phi(q) \beta \ln( p_{h-class} )$$  This  could be perhaps used to associate phase variations along $t$ of $L(s,\chi)$,  ( or  $ \left| \ln[\L(s,\chi)] \right|$ ),   peaks to particular $h-$classes of primes }
 ,  we will show that ( $\forall t \ , \forall \epsilon >0$) :
 \be \label {Fourth} 
  \left| 
 \sum_{  p<p_{max}   \    } 
 \int_{p- \alpha \ln(p) }^{p + \beta \ln(p)}  
 \left\{  
 \frac{ \cos [\ln(p) t  ]   } { (\alpha+\beta) \ln(p) p^{1/2+\epsilon}}
  -
   \frac{\cos(\ln(y) \ t    )} {y^{1/2+ \epsilon}}     \frac {1}{\ln(y)}  \right\} dy 
   \right|
    <
    \sum_{  p<p_{max}     } 
 \int_{p- \alpha \ln(p) }^{p + \beta \ln(p)}  
 \left\{  
 \frac{1 } 
 {  p^{1/2+\epsilon}    y^{1/2+ \epsilon}   }
   \right\} dy   < \infty
 \ee

\noindent So, as the computation for imaginary part is similar,  we have convergence  of expression  $\ln[ EP(s)] - \int_0^\infty x^{-s} d[ Li(x)]$ for $\epsilon>0$. 
 
 \noindent    We can write  using  (\ref   {parteImportante} ):

 \be \label {ApprMinus0_1}
% -\frac {\ln(p^*)}{ \pi} 
 \sum_{  p<p_{max}     } 
 \int_{p- \alpha \ln(p) }^{p + \beta \ln(p)}  
 \ee
 $$
 \left\{  
 \frac{ \cos [\ln(p) t  ]       } { (\alpha+\beta) \ln(p) p^{1/2+\epsilon}}
  -
   \frac{\cos(\ln(y) \ t    )} {y^{1/2+ \epsilon}}    % \sin\left(  \frac{\pi \ln(y)}{\ln(p^*)} \right)     
    \frac {1}{\ln(y)}  \right\} dy =
  $$

   \be \label {ApprMinus0_2}
% -\frac {\ln(p^*)}{ \pi} 
 \sum_{  p<p_{max}     } 
 \int_{p- \alpha \ln(p) }^{p + \beta \ln(p)}  
 \ee
 $$
 \left\{  
 \frac{
  y^{1/2+ \epsilon}  \ln(y)      \cos [\ln(p) t ]    %\sin\left( \frac{\ln(p) \pi}{\ln(p^*)}   \right)   
     -  (\alpha+\beta) \ln(p) p^{1/2+\epsilon}  \cos(\ln(y) \ t )     % \sin\left(  \frac{\pi \ln(y)}{\ln(p^*)} \right)  
  } 
 { (\alpha+\beta) \ln(p) p^{1/2+\epsilon}    y^{1/2+ \epsilon}  \ln(y) }
   \right\} dy 
  $$
  \noindent Inside the interval  $p- \alpha \ln(p) < y <p + \beta \ln(p) $ ( $ \Delta y= y-p$) ,  
  to first order we have:
  $$ \cos(\ln(y) \ t )  =    \cos(\ln(p) \ t ) -\frac{t(y-p)}{p}\sin( \ln(p) \ t )$$
 
 % $$ \sin\left(  \frac{\pi \ln(y)}{\ln(p^*)} \right) = \sin\left(  \frac{\pi \ln(p)}{\ln(p^*)} \right) + \frac{\pi (y-p)}{ p \ln(p^*)}    \cos\left(  \frac{\pi \ln(p)}{\ln(p^*)} \right)$$
 
 so { \bf at first order } in  $\frac{(y-p)}{p}$, with $p >> Max(\alpha, \beta) \times \ln(p)$:
 
 $$  y^{1/2+ \epsilon}  \ln(y)      \cos [\ln(p) t  ]   %\sin\left( \frac{\ln(p) \pi}{\ln(p^*)}   \right)   
    -  (\alpha+\beta) \ln(p) p^{1/2+\epsilon}  \cos(\ln(y) \ t  )  
    % \sin\left(  \frac{\pi \ln(y)}{\ln(p^*)} \right) 
     =
 $$
 $$
 \cos(\ln(p) \ t )   % \sin\left(  \frac{\pi \ln(p)}{\ln(p^*)} \right)
  [ y^{1/2+ \epsilon}  \ln(y) -   (\alpha+\beta) \ln(p) p^{1/2+\epsilon} ]  +
 $$
 $$
 \frac{(y-p)}{p} \left\{  % \frac{\pi }{  \ln(p^*)}    \cos\left(  \frac{\pi \ln(p)}{\ln(p^*)} \right) 
 -t\sin( \ln(p) \ t  )
  \right\}
 $$

 So the numerator of the fraction under integral in (\ref{ApprMinus0_2}) can be split in two parts resulting in two integrals:

  \be \label {FirstO} 
  \frac {\ln(p^*)}{ \pi}  
  \sum_{  p<p_{max}      } 
 \int_{p- \alpha \ln(p) }^{p + \beta \ln(p)}  
  \frac{          % \frac{\pi }{  \ln(p^*)}    \cos\left(  \frac{\pi \ln(p)}{\ln(p^*)} \right) 
 -t\sin( \ln(p) \ t   ) } 
  { (\alpha+\beta) \ln(p) p^{1/2+\epsilon}    y^{1/2+ \epsilon}     \ln(y) } \ 
    \frac{(y-p)}{p} dy   
    \ee
 
 that in absolute value is less than
 
  $$
 % \frac {\ln(p^*)}{ \pi} 
   \sum_{  p<p_{max}      } 
 \int_{p- \alpha \ln(p) }^{p + \beta \ln(p)}  
  \frac{   t% +1
  } 
  { (\alpha+\beta) \ln(p) p^{1/2+\epsilon}    y^{1/2+ \epsilon}     \ln(y) } \ 
    \frac{(\alpha+\beta) \ln(p) }{p} dy   = $$
    $$
  %\frac {\ln(p^*)}{ \pi}    
   \sum_{  p<p_{max}      } 
 \int_{p- \alpha \ln(p) }^{p + \beta \ln(p)}  
  \frac{   t %+1
   } 
  {  p^{3/2+\epsilon}    y^{1/2+ \epsilon}     \ln(y) } \ 
    dy   
     $$
 
 Wich is surely  convergent for $p_{max}\rightarrow \infty$  ,  $\forall t$  for $\epsilon>0$  (and beyond).
 
 And

  \be \label {prinpart} %\frac {\ln(p^*)}{ \pi} 
    \sum_{  p<p_{max}     } 
 \int_{p- \alpha \ln(p) }^{p + \beta \ln(p)}  
  \frac{  \cos(\ln(p) \ t   )      % \sin\left(  \frac{\pi \ln(p)}{\ln(p^*)} \right)
  [ y^{1/2+ \epsilon}  \ln(y) -   (\alpha+\beta) \ln(p) p^{1/2+\epsilon} ]  } 
  { (\alpha+\beta) \ln(p) p^{1/2+\epsilon}    y^{1/2+ \epsilon}  \ln(y) }
   dy   \ee
 
\noindent We observe that  
 $\left|    \cos(\ln(p) \ t )   %\ \frac {\ln(p^*)}{ \pi}  \ \sin\left(  \frac{\pi \ln(p)}{\ln(p^*)}\right)
  \right | \le 1$ is upper bouded by 1,
 and :
 
 $y^{1/2+ \epsilon}  \ln(y) -   (\alpha+\beta) \ln(p) p^{1/2+\epsilon} $ in the interval  is upper bounded by 
 $(\alpha+\beta) \ln(p)    \ln(y) $. 
 
 Because,  in interval  $p- \alpha \ln(p) < y <p + \beta \ln(p) $, 
 
 we have:
 
   $|y^{1/2+ \epsilon}   -    p^{1/2+\epsilon} | <\ln(y) < Max(,\alpha,\beta) \ln(p) \ \ ; \forall \epsilon:  0<\epsilon \le 1/2 $ 
 
 but 
 
  $ |  \ln(y) -   (\alpha+\beta) \ln(p) |   <  (\alpha+\beta) \ln(p) $

  so,  in interval  $p- \alpha \ln(p) < y <p + \beta \ln(p) $: 
  
  $\left| y^{1/2+ \epsilon}  \ln(y) -   (\alpha+\beta) \ln(p) p^{1/2+\epsilon} \right |  <(\alpha+\beta) \ln(p)    \ln(y) $
  
%  DA RIVEDERE !!!!!!

\noindent Orders greater than the first    (i.e. all terms but the first in (\ref {LogEP2}) )   lead to  convergent integrals  %like their sum
 ( at least for $\epsilon>0$) like (\ref{FirstO}).

 \noindent  So (\ref  {prinpart}) in absolute value  is upper  bounded by 
 
  \be \label {ApprMinus0_3}
 % \frac {\ln(p^*)}{ \pi} 
 \sum_{  p<p_{max}     } 
 \int_{p- \alpha \ln(p) }^{p + \beta \ln(p)}  
 \left\{  
 \frac{1 } 
 {  p^{1/2+\epsilon}    y^{1/2+ \epsilon}   }
   \right\} dy 
  \ee

  \noindent This means that for $\epsilon>0$ , i.e $\Re(s)>1/2$,  left part of (\ref {Third}), i.e. expression: $\ln[ EP(s)] - \int_0^\infty x^{-s} d[ Li(x)] $, is absolutely convergent. END OF PROOF
  
 %    $$ - - -$$

%    $$---from-appendix$$

%%%%%%%%%%%%%%%%%%%%%%%%%%%%%%%%%%%%%%%%                       
\subsection{ Lemma2 B: Euler product in critical strip }   \label {EPInCR}

  \noindent Let us consider the expressions:                     
                        
                        $$\prod_{\forall p} \frac{1}{1-\frac{1}{p^s}} \quad p \quad prime  \quad \ \ \; \ \ \ 
                        \ \zeta(s.)=\sum _{n=1}^\infty  \   n^{-s}
                        $$

 \noindent     It is  known that    $\prod_{\forall p} \frac{1}{1-\frac{1}{p^s}} \quad p \quad prime  \quad \ $ can be seen as the product of infinite geometrical serie, each one,  with common ratio $ \frac{1}{p_j^s}$, ( with  $ \left|   \frac{1}{p_j^s} \right| <1  $ ) i.e. $  \frac{1}{1-\frac{1}{p^s}} = \sum_{n=0}^\infty  \left(\frac{1}{p_j^s} \right)^n$. Let us consider only a finite number of primes  $j_{max}: \  p_1, p_2 . . .p_{max} \ , \ i.e.  j_{max}=\pi(p_{max}) $, and so also a finite numbers of geometrical series to be multiplied. Besides notice that the sum of  geometrical series can be written by a finite number of terms:

\noindent                     We can multiply $ j=1,...j_{max}=\pi(p_{max}) $ geometric series with infinite terms, but  we can alternatively choose  to consider, in each geometrical series,   the sum beyond  a certain  exponent $\alpha_j =\left \lceil \frac {\ln(p_{max})}{\ln( p_{j })} \right \rceil$ as a whole. i.e.
           $$ \frac{1}{1-\frac{1}{p_j^s}} = \sum_{n=0}^\infty  \left(\frac{1}{p_j^s} \right)^n=
           \sum_{n=0}^{\alpha_j-1}  \left(\frac{1}{p_j^s} \right)^n + \left(\frac{1}{p_j^s} \right)^{\alpha_j} ( 1-p_j^{-s})^{-1}
           $$         
    So we can have a finite number of terms for each of the    $ j=1...j_{max}$ geometric series instead of infinite terms.

 %We can arrange the result as follow 
 \noindent For example as an $n \approx  p_{max}$ cannot have two divisors both $ >\sqrt{p_{max}}$   then all the  series with $p_j$ from $\approx \sqrt{p_{max}}$ to $p_{max}$ are  simplified as :

   $$1+ \frac{1}    { p_{j}^s}   + \left(  \frac{ 1}    { p_{j }^s} \right)^2 \frac{1}{1-\frac{1}    { p_{j }^s}}
   \   \  ; \  \sqrt {p_{max}}   < p_j \le   p_{max}   $$
   
   \noindent At the end we get:
   $$EP(t,\epsilon,p_{max})= $$
       
                       \be \label {EulerTrickFiniteTerms}
                      \prod_{p=2}^{p_{max}}\frac{1}{1-\frac{1}{p^s}}
                       = \left\{
                       \sum_{n=1}^{p_{max}}\frac{1}{n^s} 
                       \right\}
                       +
                       \left\{
                       \sum_{\rho>p_{max}} \frac{1}{\rho^s} \prod_{\alpha_j >0} ( 1-p_j^{-s})^{-1}  
                        \right\}
                        = \zeta(s,p_{max} )+R(s,p_{max})                 
                       \ee
                       $$
                       Where  \  \rho= p_1^{\alpha_1} p_2^{\alpha_2} . . .p_{max}^{\alpha_{max}}  \ ;  \alpha_j = 0,1,2 . . .  \   \  and \ \  j_{max} = \pi(p_{max})$$

                      \noindent So we can write: 
                      \be \label {ZfromEP}    \zeta(s,p_{max} )= EP(t,\epsilon,p_{max}) -R(s,p_{max})              \ee
                                         
                    \noindent   In (\ref{EulerTrickFiniteTerms}) we can have huge amount of terms :

                     \be \label {NOfTerms} N^{o } (terms)=  3 \prod_{j=1}^{j_{max}-1}   
                        % \left [ \frac{\ln(p_{max})} {\ln(p_j)} \right]
                        \left(
                         \left  \lceil   \frac{\ln(p_{max})} {\ln(p_j)} \right  \rceil +1
                         \right)
                     %   \left \lceil   \right \rceil
                     %   \left \lfloor  \right \rfloor
                          %  \ \ - p_{max}
                        \ \ ; \ \ 
                        %  \prod_{p} 
                        last \  factor
                         \left  \lceil   \frac{\ln(p_{max})} {\ln(p_{max})} \right  \rceil +1   \ \ ; \ \  is  \ put \ to \ 3
                          \ee

                          %,  because for $p>p_{max}$ the factor is $1$.
                        %    \noindent Anyhow a finite number for finite $p_{max}$.    
          
           \noindent  For example for $p_{max}= 31$ we have in  (\ref{NOfTerms} )  that $ N^{o } (terms)= 787320$, whoose only 31 belong to $\zeta(s,p_{max} )$. The others belong to $R(s,p_{max})=R(t,\epsilon,p_{max})$. While $\rho_{max}(p_{max} )$, i.e. the greatest $\rho$  in (\ref{EulerTrickFiniteTerms}), is $>> Primorial(p_{max} )$, so an huge number too, though finite. Besides if we consider a successions of 
           \be \label {Succ} 
           p_{max}^k = \rho_{max}^{k-1}(p_{max}^{k-1})   \ \ ; where \ k \ is \ index \ not \ power
           \ee
            it is apparent  that all the terms of  the sum  in $R(t,\epsilon,p_{max}^{k-1})$, are absent in $R(t,\epsilon,p_{max}^{k})$. So it is apparent that a build-up toward a certain value, from whatever $p'_{max}$ to  $p''_{max} \rightarrow \infty $,  seems to be excluded for  $R(t,\epsilon,p''_{max})$ when  $p''_{max} \ge \rho_{max}(p'_{max})$.  And we are speaking always of huges but  finite set of terms.
          
    %  \footnote{ Lemma: if $R(t,\epsilon,p_{max}) \rightarrow \lambda <\infty $, for $p_{max} \rightarrow\infty $  then necessarily  $\lambda=0$. PROOF: for  $p_{max}$ big enough , $R(t,\epsilon,p_{max})$    is  close as we want to $\lambda$.  So the sum of the new terms from $p'_{max}$  to $p''_{max}$ can be small at will. But if  $p''_{max}=\rho_{max}(p'_{max} )$, then  the added new terms  are  exactly  $R(t,\epsilon,p''_{max})$, so we fulfill the  starting hypothesis  only if  $\lambda=0$  }     

       %   $$ - - - -$$

    % $$ - - - -15-03-26 -A$$

\noindent Which is  the limit of $R(t,\epsilon,p_{max} )\ $ for $p_{max} \rightarrow \infty$?  Considering a subsuccession like (\ref {Succ} ) one can presume also no limit at all because the terms  of (\ref {Succ} )  are completely decorrelated. In general $R(t,\epsilon,p_{max} )\ $ for $p_{max} \rightarrow \infty$ can have a limit $\lambda$ ( finite or infinite)  or no limit at all. Lukily  the term  $R(t,\epsilon,p_{max} )\ $ is embedded in formulas where, thank to (\ref {ZfromEP}), we know the behavior  of all terms but  $R(t,\epsilon,p_{max} )\ $.  % (also in critical strip, formally , without presuming convergence)

{ \bf Lemma 2B : For $p_{max}\rightarrow \infty$,  $R(s,p_{max})$ tends to zero at least for $\epsilon>0$.}

PROOF:
%  $$\zeta(s,p_{max})= EP(s,p_{max})- R(s,p_{max})$$

%and $\zeta(s,p_{max})$ enters in 

$$ \ln(\zeta(s,p_{max}) - \int_0^{p_{max}} x^{-s} d[Li(x)] 
%+ \ln \left(1-\frac{R(s,p_{max})}{EP(s,p_{max})} \right)
\ \ \ \ \ \ \ \ \ \
\rightarrow  \ \ \ \ \ \ \ \ \ \ \
\ln[\zeta(s)(s-1)] 
$$
that { \bf  is valid also on critical strip,}, can be written  as:

\be \label {EPandZ} \ln(EP(s,p_{max}) - \int_0^{p_{max}} x^{-s} d[Li(x)] + \ln \left(1-\frac{R(s,p_{max})}{EP(s,p_{max})} \right)
\rightarrow 
\ln[\zeta(s)(s-1)] 
\ee

%This  must be valid also on critical strip.
\noindent As we know that 
\begin{itemize}
\item $\ln(EP(s,p_{max}) - \int_0^{p_{max}} x^{-s} d[Li(x)] $ converges for  $\epsilon>0$ , see subsection \ref {conv2_12}
\item $\ln[\zeta(s)(s-1)] $  is surely defined  (excluding singular points) for $\epsilon>0$  
\end{itemize}
So we deduce that $\ln \left(1-\frac{R(s,p_{max})}{EP(s,p_{max})} \right)$ must converge too.
But as $EP(s,p_{max})$ does not converge in critical strip only option left  is $R(t,\epsilon,p_{max} ) \rightarrow0 $ for $p_{max} \rightarrow \infty$. END of PROOF
 
 \noindent Corollary: $\ln[\zeta(s)(s-1)] $ converges for  $\epsilon>0$

   \noindent  From Lemma2A and B,   could be proved  directly  RH without the help of  Angular Momentum  logic path. The crucial point   is  that (\ref {EPandZ} ) , i.e. $\ln[\zeta(s)(s-1)] $,  cannot diverge for $\epsilon>0$ (Lemma2B).
   
    \noindent On the other hand, to go on  with the exploration of Angular Momentum properties of $\zeta(s)$ and $\xi(s)$  functions, with a more founded mathematical base (at least for $\epsilon>0$), is an analysis that   presents interesting aspects for the comprehension of prime number distribution beyond possible RH proof approaches. %: for example the comparison of computed with forecasted correlations.

                                                                                                                                                                                                                                                                                                                                                                                                                                     \noindent So we prove  in appendix \ref {Conv3_21}, with same technique of Lemma 2A, that also  (\ref {ProVarPhaseCsi}) converges for $\epsilon=0$, and, as $R(s,p_{max}) \rightarrow 0 $ for $p_{max} \rightarrow \infty$, then, (\ref {ProVarPhaseCsi})  tends to  $\left\{  \frac{ \partial \Im [\ln \{ \zeta(s)(s-1) \}]}{\partial t} \right\} $ . Also this results get  directly to RH because  a zero in the upper half of critical strip would lead to a diverging $\left\{  \frac{ \partial \Im [\ln \{ \zeta(s)(s-1) \}]}{\partial t} \right\} $.

                                                                                                                                                                                                                                                                                                                                                                                                                                     \noindent So we  will use  appendix \ref {Conv3_21}  as a mathematical justification  to follow the old proposed logical path of Angular Momentum also if   and Lemma2A and Lemma 2B would lead us directly to RH.

%%%%%%%%%%%%%%%%%%%%%%%%%%%%%%%%%%%%%%%%%%%%%%%
%%%%%%%%%%%%%%%%%%%%%%%%%%
\section {Consequences of Lemma 2 of subsection \ref{CONJ}}

%\noindent  In the whole section it is assumed the conjecture of section \ref{CONJ}.

%%%%%%%%%%%%%%%%%%%%%%%%%%%%%%%%%%%%%%%%%   

 \subsection                         { LEMMA \ 3 : \ $ \frac {\partial \angle[\xi(t,\epsilon )]  } {\partial t} $ and Euler product    }  \label {Theorem1} \

 \noindent   As  (\ref  {ProVarPhaseCsi}),and   (\ref{FromJminusLi}),  are computed as mean value in the interval : $t\pm \frac{\pi}{\ln(p^*)}$  let us take all $t$ in intervals  not intersecting the intervals:  $t_k^* \pm \frac{\pi}{\ln(p^*)} \ , \  \forall k$.  Where $ \zeta( t^*_k,\epsilon^*_k)=0$) , i.e. $|t -t^*_k| > \frac{2 \pi}{\ln(p^*)} $, avoiding so the phase discontinuity  of logarithm  % 
 $ \forall \epsilon:  -|\delta| < \epsilon < \epsilon^*_k$. By  (\ref  {DeFaseSuL1SuDeT})  so, neglecting $O(t^{-2})$, and, using (\ref{FromJminusLi}), or,  (\ref  {ProVarPhaseCsi}): 
    
  \be \label {Uguale0Bis} 
  \ln\left[ \sqrt{\frac{t }{2 \pi}} \right]+\left\{  \frac{ \Delta \Im [\ln \{ \zeta(s)(s-1) \}]}{\Delta t} \right\}_{p_{max}>>p^* \rightarrow \infty}
   \rightarrow
  \frac {\partial \angle[\xi(t,\epsilon )]  } {\partial t} 
   \quad ; 
   |t -t^*_k| > \frac{2 \pi}{\ln(p^*)}   \ , \ \quad if\quad \zeta( t^*_k,\epsilon^*_k)=0
  \ee

   \noindent  While for   $ \forall \epsilon > \epsilon^*_k$ we can drop condition $ |t -t^*_k| > \frac{2 \pi}{\ln(p^*)}$. %.
  For $\epsilon=0$ \footnote{Of course far from diverging points  that are the zeros of $\zeta(s)$}, expression  (\ref{Uguale0Bis}), far from zeros, tends to zero. See subsection \ref {ZeroOnCritLine}.
 \noindent  If we increase $p^*$ %
and  $p_{max}$ (\ref{FromJminusLi}), i.e. (\ref  {ProVarPhaseCsi}) %( , 
  converges in     $(  t^*_k +\frac{2 \pi}{\ln(p^*)}   , t^*_{k+1}  -\frac{2 \pi}{\ln(p^*)}  ) \ \forall k$.

\

  \noindent  Besides in order to counter  the fictitious oscillation of Gibbs phenomenon,  stemming from the inevitable $p_{max}  <\infty$ condition (see  \cite [p.~25]{Pinsky2002}),   we could apply a moving window  of  length  $\Delta t =W= \frac{2 \pi}{\ln(p_{max})}  $. So the effective convergence  interval of  (\ref {Uguale0Bis} ) becomes:

   $$\left(  t^*_k +\frac{2 \pi}{\ln(p^*)} +\frac{2 \pi}{\ln(p_{max})}  \ \   , \ \   t^*_{k+1}  -\frac{2 \pi}{\ln(p^*)}-\frac{2 \pi}{\ln(p_{max})}   \right)$$
  
  \noindent   For  $p^*$ and $p_{max}$  big enough  we find  open interval  $(  t^*_k  , t^*_{k+1}  )$.  See fig. \ref  {SmthT350T360Pstar6Pmax60}.  %
%

   %%%%%%%%%%%%%%%%%%%%%%%%%%%%%%%%%%%%%%%%%%%
  %%%%%%%%%%%%%%

%

\subsection { LEMMA  4: estimate of sign of  $\frac{\partial^2 \Im [ \ln(\zeta(s)(s-1)) ]}{\partial \epsilon  \partial t}$ } \label {Lemma4}

\noindent %{ 
 Suppose $\Delta \epsilon >0$ but arbitrarily small, $ |t -t^*_k| > \frac{2 \pi}{\ln(p^*)}$, $ p_{max} \ge p^* \rightarrow \infty $.
   \noindent Then, for (\ref {Uguale0Bis} ) we have that (\ref  {ProVarPhaseCsi}) 
 evaluated in $\epsilon=0$  \footnote {we suppose to be far from $\zeta(s)$  zeros} is tending to $ - \ln\left[ \sqrt{\frac{t }{2 \pi}} \right] <0$ as  $ p_{max} \ge p^* \rightarrow \infty $. We affirm that the same  (\ref  {ProVarPhaseCsi}), or (\ref {exprToLookAt_2}),  evaluated in $\epsilon=\Delta \epsilon >0$  is always negative but with lower absolute value, so that  (\ref {Uguale0Bis} )  has a positive increment for  $\Delta \epsilon >0$. This is even more true  if we consider (\ref{DeFaseSuL1SuDeT}), and, (\ref {Uguale0Bis}) with $O(t^{-2})$ term:  i.e. $+\frac{3 \Delta \epsilon}{4 t^2}$.

  \noindent FOREWORD of the proof.                                                    
     
\noindent Notice that  if we use  % the approximated   (\ref {appr}) to 
(\ref  {ProVarPhaseCsi})  % %
we are in  same situation 
like in  \cite  [p.~10] {Giovanni Lodone Nov2024} if we pose the congruence modulus  $q=1$, so we could use exactly the same proof. The difference is only that in  \cite  [p.~10] {Giovanni Lodone Nov2024} the logarithmic integral   is used only to allow the introduction of PNT , while here, similar integral is crucial to allow also conditional convergence. 

\noindent We will use here   instead    (\ref{FromJminusLi}), that is exact in the sense that equi-sign intervals in  $p$ space (i.e. (\ref  {HalfRotation2_} ) and (\ref  {HalfRotation_}  )  )  are the same  both for   prime powers and for $Li(x)$ integral.. The (\ref {appr}) instead is useful in order to shed light on the primes spectrum, i.e.  which are the ``resonance'' frequencies  $\frac{\ln(p)} {2 \pi}$ in prime distribution. See fig. \ref {T710Eps01}. % apply Lemma 2 to  (\ref  {ProVarPhaseCsi}) that   is the same as   (\ref{FromJminusLi}) and both tend to $\frac{\partial \Im [ \ln(\zeta(s)(s-1)) ]}{ \partial t}$ %conjecturally at least from $\epsilon > -|\delta|$. See section \ref  {CONJ} .  

 \noindent Lemma 3  PROOF.

\noindent Let us take a close view at $\cos (\ln(p) \ t)$ oscillations in  (\ref{FromJminusLi}). %
The expression is:
       \be \label {exprToLookAt}
   \frac{\ln(p^*)} { \pi}\int_2^{p_{max}} 
     \frac{\cos(\ln(x) \ t)}{x^{1/2+ \epsilon}} \sin\left(  \frac{\pi \ln(x)}{\ln(p^*)} \right)  
   d \left [  J(x)-Li(x) \right]   %
   \quad   \ , \ with \  \ \ln(p^*) \ big
\ee                                                         
                                                                   
  \noindent Supposing, at first,      $\sin\left(  \frac{\pi \ln(x)}{\ln(p^*)} \right) >0$., i.e. $[p^*]^m                                                                  <p < [p^*]^{m  +1}$ with $m=0,2,4 . .$ (see  \cite  [p.~10-11]  {Giovanni Lodone Nov2024} )

                                                                                                 \begin{itemize}

\item Zero transition at increasing $\cos[\ln(x)t]$ : %$

%%%   {TransizMenoPiMezzi__}  at increasing  cosine,     {TransizPiuPiMezzi__}  decreasing  cosine

   \be \label {TransizMenoPiMezzi__} 
     \quad x_{0t}(k)= e^{ (2 \pi k- \pi/2) / t}   \quad increasing \quad cosine
   \ee
\item Zero transition  at decreasing   $\cos[\ln(x)t]$ : %
    
    \be \label {TransizPiuPiMezzi__} 
     \quad x'_{0t}(k) = e^{ (2 \pi k+\pi/2) / t}   \quad decreasing \quad cosine
     \ee
        \end{itemize}   

 Positive  contribution for (\ref{exprToLookAt}) integrand , %
 is given by prime powers inside  the exponentially growing (with $k$ ) interval  $\Delta x_k$ (subscript $0t$ means `zero transition'):
                                                                                                                                                                                                                                                                                                                                                                                           
 \be \label {HalfRotation_}
 x'_{0t}(k) -  x_{0t}(k)=    x_{0t}(k)  [e^{ ( \pi  / t)}-1]\approx  \frac {\pi}{t} e^{ (2 \pi k -\pi/2) / t} =  \frac {\pi}{t} x_{mean}=\Delta x_k
 \ee
       while, negative contribution  %
       is given by prime powers inside the exponentially growing (with $k$) interval $\Delta x'_k$:

      \be \label {HalfRotation2_}
   x_{0t}(k+1)-x'_{0t}(k)=   x'_{0t}(k) [e^{ ( \pi  / t)}-1]
      \approx \frac {\pi}{t} e^{ (2 \pi k+\pi/2) / t} =  \frac {\pi}{t} x'_{mean}=\Delta x'_k
   \ee

\noindent So each $t$ choice  partitions prime powers in two sets.  The positive, and, the negative contributors to  (\ref  {exprToLookAt}).

Looking at  (\ref{exprToLookAt})  we can split the integral in subintervals: $ . . \Delta x_{k-1}, \Delta x'_{k-1}, \Delta x_{k}, \Delta x'_{k}, \Delta x_{k+1}, \Delta x'_{k+1}, . . .$.  See (\ref {HalfRotation2_}), and  (\ref {HalfRotation_}).

\noindent    In   $\Delta x_{k}$ intervals ( (\ref {HalfRotation_}))  $J(x)$  gives a positive contribution: 
\noindent

 \be \label {OJK+}  
   O^+_{J} (\epsilon)_ {k}= \left|
  \frac {\ln(p^*)} { \pi}  
           \int_{ p >x_{0t} (k)}^{ p <x'_{0t} (k)} 
            \frac{ 
            %\ 
              \cos [\ln(p)t]\sin [\pi \ln(p) /\ln(p^*)] }{\sqrt{p} p^\epsilon}   dJ(p)
              \right|
              \ee
 
 \noindent while $Li(x)$ a negative one:
 
\noindent 
 \be \label {OLiK-}  
   O^-_{Li} (\epsilon)_ {k}= \left|
  \frac {\ln(p^*)} { \pi}  
           \int_{ p >x_{0t} (k)}^{ p <x'_{0t} (k)} 
            \frac{ 
            %\ 
              \cos [\ln(p)t]\sin [\pi \ln(p) /\ln(p^*)] }{\sqrt{p} p^\epsilon}   d[Li(p)]
              \right|
              \ee

\noindent  In   $\Delta x'_{k}$ intervals $J(x)$ gives a negative contribution:

 \be \label {OJK-}  
   O^-_{J} (\epsilon)_ {k}= \left|
  \frac {\ln(p^*)} { \pi}  
           \int_{ p >x'_{0t} (k)}^{ p <x_{0t} (k+1)} 
            \frac{ 
            %\ 
              \cos [\ln(p)t]\sin [\pi \ln(p) /\ln(p^*)] }{\sqrt{p} p^\epsilon}   dJ(p)
              \right|
              \ee
 
\noindent   while $Li(x)$ a postive contribution:
 \be \label {OLiK+}  
   O^+_{Li} (\epsilon)_ {k}= \left|
  \frac {\ln(p^*)} { \pi}  
           \int_{ p >x'_{0t} (k)}^{ p <x_{0t} (k+1)} 
            \frac{ 
            %\ 
              \cos [\ln(p)t]\sin [\pi \ln(p) /\ln(p^*)] }{\sqrt{p} p^\epsilon}   d[Li(p)]
              \right|
              \ee

\noindent The four above are positive definite quantities. Signs are made  explicit when they are assembled together.  Then  (\ref{exprToLookAt}), using (\ref {OJK+}), (\ref  {OLiK-} ), (\ref {OJK-} ), (\ref   {OLiK+}) , 
can be written also as :

     \be \label {exprToLookAt_2}
   \frac{\ln(p^*)} { \pi}\int_2^{p_{max}=x_{k_{max}+1}} 
     \frac{\cos(\ln(x) \ t)}{x^{1/2+ \epsilon}} \sin\left(  \frac{\pi \ln(x)}{\ln(p^*)} \right)  
   d \left [  J(x)-Li(x) \right]   =
\ee

  $$\sum_k^{k_{max}}\{   O^+_{J} (\epsilon)_ {k}  -    O^-_{Li} (\epsilon)_ {k}  - O^-_{J} (\epsilon)_ {k}  + O^+_{Li} (\epsilon)_ {k} \}$$

\noindent  {\bf Grouping must not alter the  order of the sum because there is not absolute convergence}. For our choices, in 
  the $k-$th oscillation, in the space of primes, %
    $ O^+_{J} (\epsilon)_ {k} $ comes before %
    $ O^-_{J} (\epsilon)_ {k} $, and , $O^-_{Li} (\epsilon)_ {k}$ comes before $O^+_{Li} (\epsilon)_ {k}$.
 \noindent  For Heath-Brown condition, \cite  {Heath-Brown : 1988}  , and   \cite  [p.~12]  {Giovanni Lodone Nov2024} , and $k$ big
 ,  we have 
 for example :

 \be \label {Rapporti}    \left[ \frac{ O^+_{J} (\epsilon)_ {k} }{  O^-_{Li} (\epsilon)_ {k}} \right]_{k \ big} \rightarrow 1
   \ \ \ or \ \ \ 
     \left[ \frac{O^-_{J} (\epsilon)_ {k} }{ O^+_{Li} (\epsilon)_ {k}} \right]_{k \ big} \rightarrow 1
   \ee.

    \noindent  Referring   to (\ref {Rapporti})   let us take the following ratio:%%

   \be \label {RappSomme}\frac{\sum_+^{p_{max}} (\epsilon)}{ \sum_-^{p_{max}} (\epsilon) }=
    \frac
    {\sum_k O^+_{J} (\epsilon)_ {k}  + O^+_{Li} (\epsilon)_ {k} }
    {\sum_k  O^-_{Li} (\epsilon)_ {k} +O^-_{J} (\epsilon)_ {k}   } \ \ \  \rightarrow  \ \ 1
   \ee
   
\noindent Where $ \sum_\pm^{p_{max}} (\epsilon)$ are the total positive or negative contribution to (\ref {exprToLookAt_2}) till $p_{max}$ , at given $\epsilon$ ( taking in consideration the sign of   $\sin\left(  \frac{\pi \ln(x)}{\ln(p^*)} \right) $, i.e. the parity of $m$ in $[p^*]^m                                                                  <p < [p^*]^{m  +1}$, see   \cite  [p.~13]  {Giovanni Lodone Nov2024}  ). 
   
  \noindent Of course it goes to 1 because the difference between numerator and denominator is %  
  of lower order with respect to them. {\bf This   is true in a continuous way  also considering fractions of oscillations}. Apart  the absolute value of cosine shape  , the common factor $\left| \frac{\sin( \pi \ln(p) /\ln( p^*)) }{p^{1/2+\epsilon}} \right|$ and the alternating choice  between  $ J(x)$ and $ Li(x)$. It is     similar to   P.N.T. ratio.

 \noindent  In other words  we can write P.N.T. as: 
  
  $$
   \frac
   {  \sum_k   \{    \Delta [J(x)] _{x_k} ^{x'_{k} }  +\Delta [J(x)] _{x'_k} ^{x_{k+1} }     \}    }
 {\sum_k \{ \Delta [Li(x)] _{x_k} ^{x'_{k} }  +\Delta [Li(x)] _{x'_k} ^{x_{k+1} }  \}}  \sim \frac{\pi(x_{k_{max}} )}{\Li(x_{k_{max}} )} \rightarrow 1 \ \ if \ \ k \rightarrow \infty$$

 \noindent but, alternating   between  $ J(x)$ and $ Li(x)$,  we have also:

  $$
   \frac
   {  \sum_k   \{    \Delta [J(x)] _{x_k} ^{x'_{k} }  +\Delta [Li(x)] _{x'_k} ^{x_{k+1} }     \}    }
 {\sum_k \{ \Delta [Li(x)] _{x_k} ^{x'_{k} }  +\Delta [J(x)] _{x'_k} ^{x_{k+1} }  \}}  \rightarrow 1 \ \ if \ \ k \rightarrow \infty$$

 The same cosine shape $\left|\cos(\ln(p) t)\frac{\sin( \pi \ln(p) /\ln( p^*)) }{p^{1/2+\epsilon}} \right|$ at numerator or denominator, or, ending with a fractional interval does not change the limit.

\noindent    If we change  $\epsilon' $ with $ \epsilon'' > \epsilon'  $, and we limit to a big $p_{max}<\infty$, all the positive contributions   of (\ref {exprToLookAt_2}) will be multiplied by:   
    
     \be \label {RappSommePos}  \frac
      {\sum_k O^+_{J} (\epsilon'')_ {k}  + O^+_{Li} (\epsilon'')_ {k}  }
     {\sum_k O^+_{J} (\epsilon')_ {k}  + O^+_{Li} (\epsilon')_ {k}  } = \rho^+(\epsilon',\epsilon'') \ \ \  <  \ \ 1
   \ee

    \noindent while the negative part by:

    \be \label {RappSommeNeg} \frac
     {\sum_k  O^-_{Li} (\epsilon'')_ {k} +O^-_{J} (\epsilon'')_ {k}  }
     {\sum_k  O^-_{Li} (\epsilon')_ {k} +O^-_{J} (\epsilon')_ {k}  }
    =\rho^-(\epsilon',\epsilon'')  \ \ \  <  \ \ 1
   \ee    
   \noindent  But  (\ref  {RappSommePos}) and (\ref  {RappSommeNeg})   are exactly the same 
   for  $k$ big (see (\ref  {Rapporti}) ), and,  for $\epsilon'' $ greater, but,  close enough to $\epsilon'=0$ we have that $\rho^\pm(0,\epsilon'')$ can reach 1 from below. From this point on the arguments are identical to those used in  \cite  [p.~14]  {Giovanni Lodone Nov2024}). %

  \noindent  Then we conclude:

  \begin{comment}

   \begin{itemize}
   \item the positive total sum  % 
   ( $\sum_k O^+_{J} (\epsilon)_ {k}  + O^+_{Li} (\epsilon)_ {k} $ , in the transformation %
   $\epsilon' =0 \rightarrow \epsilon''= \Delta \epsilon$ is multiplied by  $\rho^+(\epsilon',\epsilon'') $ i.e. (\ref {RappSommePos})

   \item  the negative total sum %
    ( $\sum_k  O^-_{Li} (\epsilon)_ {k} +O^-_{J} (\epsilon)_ {k} $ , in the transformation  $\epsilon' =0 \rightarrow \epsilon''= \Delta \epsilon$  is multiplied by $\rho^-(\epsilon',\epsilon'') $ i.e.  (\ref {RappSommeNeg})
    
    \item but $\rho^+(\epsilon',\epsilon'') =\rho^-(\epsilon',\epsilon'')= \rho(\epsilon',\epsilon'')  $  for $p^*$ fixed, but big, and $p_{max } \rightarrow \infty$
   
  \item  $\left\{  \frac{ \Delta \Im [\ln \{ \zeta(s)(s-1) \}]}{\Delta t} \right\}_{s=1/2+it \  , p^* , p_{max}} $, in the transformation  $\epsilon' =0 \rightarrow \epsilon''= \Delta \epsilon$ is multiplied by same factor. So we have :

 %  

  \be \label {VariazPos}
  \left\{  \frac{ \Delta \Im [\ln \{ \zeta(s)(s-1) \}]}{\Delta t} \right\}_{s=1/2+\Delta \epsilon+it \  , p^* , p_{max} \rightarrow \infty}
  =  \rho(0,\Delta \epsilon ) \times
  \left\{  \frac{ \Delta \Im [\ln \{ \zeta(s)(s-1) \}]}{\Delta t} \right\}_{s=1/2+it \  , p^* , p_{max} \rightarrow \infty}  
  \ee
  
\noindent   with  $\rho(0,\Delta \epsilon) <1$  so

  \end{itemize}

   \end{comment}

  \noindent    Lemma 4  thesis is proved.

$$ COROLLARY \ of \ Lemma \ 4 $$
 
   \noindent As Lemma 4  is valid for $p_{max}\ge p^* \rightarrow \infty$ we conclude (reasons are detailed in\cite  [p.~15]  {Giovanni Lodone Nov2024}):

  \be \label {ZSmen1}
\left\{    \frac{\partial^2 \angle[\zeta(s)(s-1)]}{\partial \epsilon \partial t} \right \}_{\epsilon=0} \ge 0  \ \ ;  \ \ \forall t \ne t^*_k  \ \ where\ \  s^*_k = \frac{1}{2}+\epsilon_k^*+ i t_k^* \ \ with \ \ \xi(s_k^*)=0
  \ee

 \subsection                        {  THEOREM    1 : variation of $\left\{ \mathcal{L} [ \xi (s) ] \right\}_{\epsilon \approx 0}$ with $\epsilon$ } \label {Theorem2}

 $$\forall t \ \ with \ \   \xi(t) \ne 0 \rightarrow \left(  \frac{d \xi(t)}{dt} \right)^2- \xi(t) \times \frac{d^2 \xi(t)}{dt^2} >0 \ \ where  \ \ \xi(t)= \xi(1/2+it) $$. 

PROOF.

\noindent If $t \ne t^*_k  \ \  \forall k$, for  Corollary 1, i.e. (\ref {ZSmen1}),
      \noindent % 
  and (\ref {DeFaseSuL1SuDeT}) : 
    $$\left\{ 
     \frac {\partial}{\partial \epsilon} \frac {\partial}{\partial t}  \angle [\xi(s)] \right\}_{\epsilon=0}
    =
    \frac{\partial}{\partial \epsilon}\left\{   \frac{\partial \angle[\zeta(s)(s-1)]}{ \partial t} +
     \frac{\partial
\angle\left[ \Gamma\left(\frac{s}{2}+1\right)\pi^{- s/2} \right] }  {\partial t}  \right\}_{\epsilon=0}
     >\frac{3 }{4 t^2}>0$$. %
      
       \noindent  From (\ref  {AngMomDef}) 
     we have : % 
     $
 \mathcal{L} [ \xi (s) ]= |\xi(s)|^2 \frac {\partial}{\partial t}  \angle [\xi(s)]$. %
  Then,    
\noindent   considering that:
    $$\left\{ \frac {\partial  |\xi(s)|^2}{\partial \epsilon} \ \frac {\partial}{\partial t}  \angle [\xi(s)]\right\}_{\epsilon=0}=0$$
    we have:
 \be \label {ProBombieriP6}
 \left\{  \frac {\partial}{\partial \epsilon}
 \mathcal{L} [ \xi (s) ] \right\}_{\epsilon=0}
 =  \left\{ 
 |\xi(s)|^2  \frac {\partial}{\partial \epsilon}\frac {\partial}{\partial t}  \angle [\xi(s)]  \right\}_{\epsilon=0} 
 >\delta  \ge |\xi(t,\epsilon = 0)|^2 \  \frac{3 }{4(t )^2} >0 , if \ \xi(t,\epsilon = 0) \ne 0
 \ee

\noindent   Suppose  now    $\xi(t^*_k,\epsilon^*_k>0) = 0$ , so  $\xi(t^*_k,\epsilon=0) \ne 0$.

 \noindent If   we  shrinks $|t-t^*_k| <\frac{ 2 \pi}{\ln(p^*)} $ interval with $p_{max}\ge p^*\rightarrow \infty$,  that (\ref  {ProBombieriP6}) is valid everywhere around $t^*_k$ for $\epsilon=0$. In other words with $p^*$ big enough the validity of (\ref{Uguale0Bis} ) can be reached close to the phase discontinuity line ($(t=t^*  \ , \  \epsilon <\epsilon^*)$)  as we want. Besides, computing with (\ref{ZSinhECosh1}),    $[\xi'(t)]^2 -\xi(t)  \xi''(t) =\left\{  \frac {\partial}{\partial \epsilon}
 \mathcal{L} [ \xi (s) ]\right\}_{\epsilon=0} $
 is continuous. 
 So (\ref  {ProBombieriP6}) follows also for $t=t_k^*$.  So thesis is proved  because, for Lemma 2,  $[\xi'(t)]^2 -\xi(t)  \xi''(t) 
 = \left\{  \frac {\partial}{\partial \epsilon}
 \mathcal{L} [ \xi (s) ]\right\}_{\epsilon=0}$ and here we have proved that 
 $\left\{  \frac {\partial}{\partial \epsilon}
 \mathcal{L} [ \xi (s) ]\right\}_{\epsilon=0} >0 $  where  $\xi(t)  \ne 0$  . %.       

 \noindent  In any case if  $\xi(t^*_k,\epsilon^*_k=0) = 0$ we have at least :
   $[\xi'(t^*_k)]^2 -\xi(t^*_k)  \xi''(t^*_k) 
 = \left\{  \frac {\partial}{\partial \epsilon}
 \mathcal{L} [ \xi (s) ]\right\}_{\epsilon=0} \ge 0 $  in an open interval of  $t^*_k $.

%%%%%%%%%%%%%%%%%%%%%%%%%%%%%%%%%%%%%%%%%%%%%%%
%%%%%%%%%%%%%%%%%%%%%%%%%%
\section {Conclusions}

Reality cannot be seen  as a single logic thread,   but it is actually a woven fabric of logic threads. Here for example:
the thread of  $\ln| \zeta(s)(s-1)|$, the thread of angular momentum, the thread of $\frac{\partial \Im[\ln(\zeta(s)(s-1))]}{\partial t}$, and, perhaps a lot of others logical paths are likely to lead to RH.

We found an expression  (\ref{FromJminusLi}) that   is a correlation between: (1) difference of prime powers distribution, $J(x)$,  from the mean $Li(x)$, and,  (2) a weighted $\cos(\ln(x) t)$. It can be interpreted  as the ``spectrum of primes'' \cite {Mazur:2015}.
  \noindent In  (\ref  {FromJminusLi})   %, 
  formulation,  %
  the zeros of $\zeta(s)$ are the $t-$abscissa of
  correlation peaks  of $d[J(x)-Li(x)]$,  %
   with the weighted  $\cos(\ln(x)t)$, along $x$.

    \noindent From a simplistic point of view one can visualize a point moving along a fixed  $\epsilon$ line toward growing $t'$ and compute the phase variations associated at each zero of $\zeta(s)$. following $\Im[\ln(\zeta(s)]$ , i.e.  (\ref {Trasf1_}). This computation can be done on imaginary part of $\ln(\zeta)$, or on average, using prime distribution,  (for $t$ far from $t'$), by (\ref{Trasf1ImPartDerivSuT}). Their difference, so, will enhance phase variation of close zero only, because, the effect of far zeros can be reasonally computed also by mean prime distribution. So  fig. \ref {T710Eps01} (or similar) are  intuitive and highly predictable.

   \noindent We are free to choose resolution  (changing $\epsilon$ and $p^*$) at will, and maximum number of primes involved ($p_{MAX}$). With $\epsilon-\epsilon^*_k \le 0$  (see (\ref {Uguale0Bis} ) ) we have that the peak at $\epsilon=\epsilon^*_k  \  , \  t=t^*_k$  becomes a ridge in  the  $ |t -t^*_k| < \frac{2 \pi}{\ln(p^*)}$ region, like in $\frac{\partial \Im[\ln(\zeta(s)(s-1))]}{\partial t}$. This ridge    shrinks, along $t$, as $p^*  \rightarrow \infty$ even if it is  prohibitive from a computational point of view. See fig. \ref  {T1640Luglio30-2}, and, fig  \ref {SmthT350T360Pstar6Pmax60}. 

\noindent From Lemma2  analysis of Euler product, in order  to strengthen  his  use  in critical strip, we encounter also arguments to conclude that RH must to hold also  directly from Euler product without  using  equivalence  \cite[p.~6] {Bombieri:2000cz}, reported here p. 2, and 
(\ref {StrongVersion} ) and angular momentum concept. But ,the  angular momentum way, leeding naturally to spectra concept, seems  to ease the comprehension of primes behavior. 

% \noindent Besides with  Theorem 1 we connect this theory with equivalence   in section (\ref{Theorem2})   stated in \cite[p.~6] {Bombieri:2000cz}.  See page 2 here.  
 
 \noindent Statement (\ref{Theorem2}) does not exclude  double zeros on critical line, as it would be if we could erase  $\xi(t)\ne0$ condition in (\ref{Theorem2}).

 \noindent The similarity with theorem in  \cite  [p.~15] {Giovanni Lodone Nov2024}  about Dirichlet odd primitive characters  is obvious, and it extends also to the short way that arises in Lemma2.
 
 \noindent This study , needless to say, is a preliminary proposal, so we'll be glad to receive whatever observation from the unlikely interested reader.

$$Acknowledgments
$$
I thank Paolo Lodone for extremely useful discussions, and, for contributing in some point of this work. Although the author bears the whole responsibility for claimed results, it is certain that, without our periodical, sometime ``fierce'', brainstormings, this paper would never come to light.
 %

%%%%%%%%%%%%%%%%%%%%%%%%%%%%%%%%%%%%%%%%%%%%%%%%%%%%%%%%%%
%%%%%%%%%%%%%%%%%%%%%%%%%%%%%%%%%%%%%%%%%
%   \section{Bibliography}

\vspace{0.3cm}

 %         

%         

%%%%%%%%%%%%%%%%%%%%%%%%%%%
%%%%%%%%%%%%%%%%%%%%%%%%
\appendix

 \section                                                                                                                                                                                                                                                                                                                                                                                                                                                                                                                                                                                                             {
Numerical  comparisons between $\frac{\partial  \angle [\xi(t,\epsilon)]}{\partial t} $ and  (\ref   {ProVarPhaseCsi})   }  

\noindent To  compute         (\ref{ProVarPhaseCsi}   )  we have extracted, by an Eratosthenes   sieve, prime numbers till $p_{max}$ i.e. about a total of  $\pi(p_{max}) \approx \frac{p_{max}}{\ln(p_{max})} $ prime numbers.  
\noindent While  (\ref {Trasf1ImPartDerivSuT}) is computed  with (\ref{IntConSinEExp} ) and  (\ref {IntConCosEExp}).
 % following indefinite integrals :
 % $ \int e^{av} \sin(bv) dv = \frac{e^{av} [a \sin(bv)-b \cos(bv) ] }{a^2 + b^2} \ \ ; \ \  \int e^{av} \cos(bv) dv = \frac{e^{av} [a \cos(bv)+ b \sin(bv) ]}{a^2 + b^2} $.

\begin{figure}[]
\begin{center}
\includegraphics[width=1.0\textwidth]{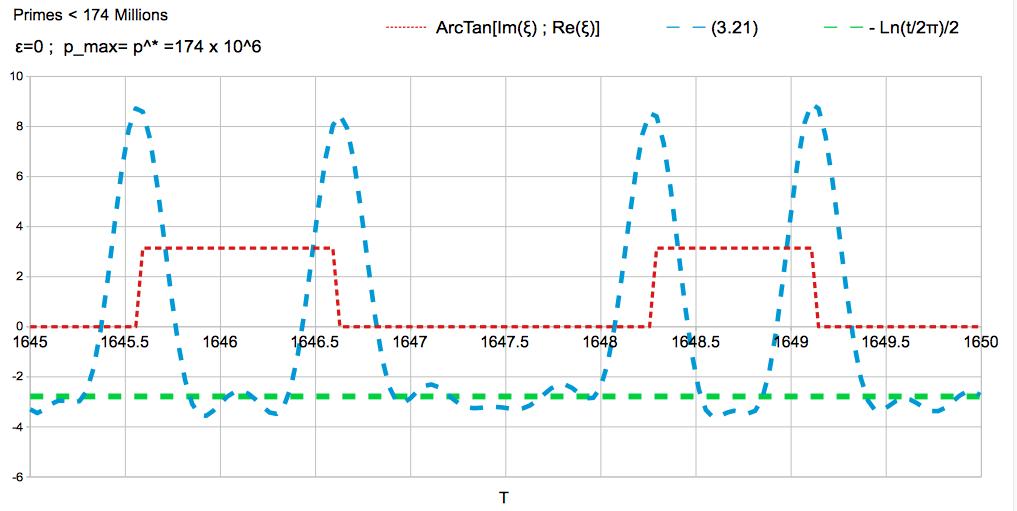}
\caption{ \small {\it  % 
Here same situation of  fig. \ref {DerAnglesVsTParamEps}  with (\ref {ProVarPhaseCsi}) at $\epsilon=0$ and  $p^*=p_{max}= 174 \times10^6$, with $\angle[\xi(t)]$ in $ 0 ,\pi$ interval.  %
Here $\Delta t=t_2-t_1=\frac {2 \pi } { \ln(174\times 10^6))}\approx 0.33$.  See (\ref{Estremi}). %
This gives an idea of the huge amounts of primes to be processed to get a closer correspondence for $\epsilon=0$. For $\epsilon >0$  less primes are needed. %
%ODS 30 07 21  DelphiEulProApproxPhaseVar
}}
\label {T1640Luglio30-2}
\end{center}
\end{figure}

\noindent  Numerical tests, from fig.  \ref  {T1640Luglio30-2} to fig. \ref {SmthT350T360Pstar6Pmax60}  ,     are intended to verify that definition of ``spectrum'' of the $J(x)-Li(x)$ function following (\ref {ProVarPhaseCsi}) %and (\ref {2Integrals_} 
presents peaks in correspondence of non trivial zeros of $\xi(s)$, and behaves  like $\frac { \partial \Im[\ln(    \xi(s)(s-1))  ]}{\partial t}  $  also with  singularity at $\epsilon < \epsilon^*_k$, that,  instead, does not appear in % . 
$\frac{\partial \angle [Z(t,\epsilon)]}{\partial t}=\frac{\partial \angle [\xi(s)]}{\partial t}$  defined by  (\ref {ZSinhECosh1}). See (\ref  {Uguale0Bis} ). In figures  \ref {T60T90xk0} and \ref  {T30T90Per3Eps} in visible the trend  toward  $- \ln\left( \sqrt{\frac{t}{2 \pi}}  \right)$   of  (\ref {ProVarPhaseCsi})  for $\epsilon=0$.

\begin{figure}[!htbp]
\begin{center}
\includegraphics[width=1.0\textwidth]{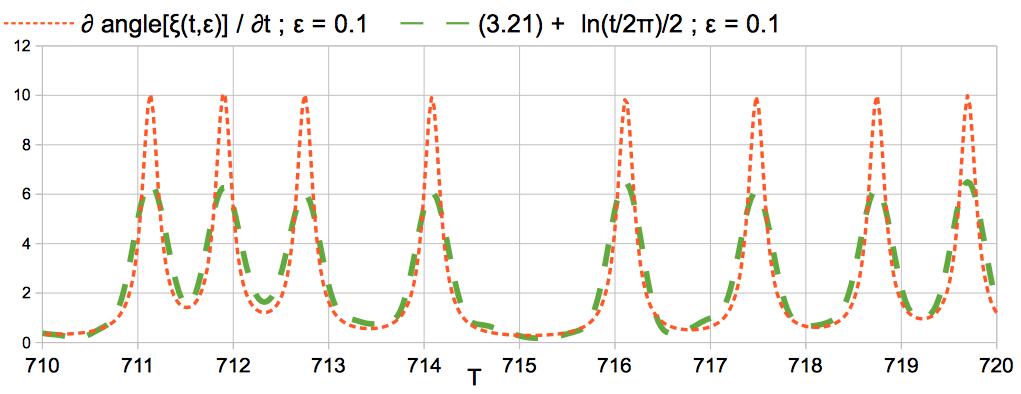} 
\caption{\small {\it   
% ODS 18_06_22ForStabilizationArgument.   sheet `` DA ODS 15 12 21''
%
  Comparison between phase variation of $\xi(t,\epsilon)$,  computed with (\ref {ZSinhECosh1}), versus  (\ref  {ProVarPhaseCsi}) summed with $ +\ln\left[ \sqrt{\frac{t }{2 \pi}} \right]$.   The value  $\epsilon=0.1$, and $p_{max} = p^* = 3 \times 10^{7}$  .  %
    When $\epsilon$ is greater than zero like in this case,  $\epsilon=+0.1$, we need  comparatively less primes, with respect to $\epsilon=0$, to attain the shape  of  $\frac{\partial \angle[\xi(s)]}{\partial t}$ computed by (\ref {ZSinhECosh1}). It is like a filtering of  high frequency ,$ \frac{\ln(p)}{2 \pi}$, i.e. high $p$ contributions. %
} }
\label {Stabiliz1}  
\end{center}
\end{figure}

\begin{figure}[]
\begin{center}
\includegraphics[width=1.0\textwidth]{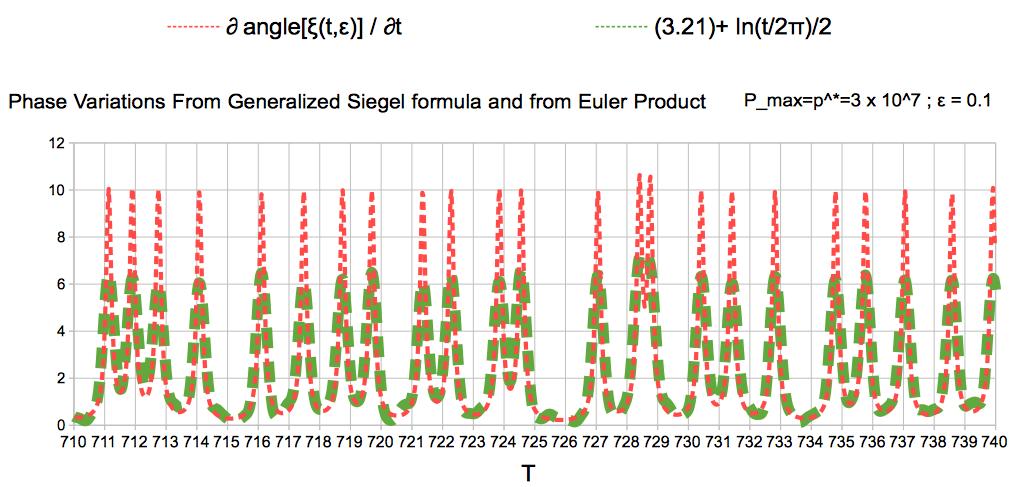}  
\caption{ \small {\it   
% ODS: 15_12_21VerificheDeduzioni
Same as fig. \ref {Stabiliz1}   with an expanded T span.
    Notice that  this $p^* $ value does not allow to resolve  two correlation spikes between $t= 728$, and, $t=729$.
}}
\label {T710Eps01}
\end{center}
\end{figure}

\begin{figure}[!htbp]
\begin{center}
\includegraphics[width=1.0\textwidth]{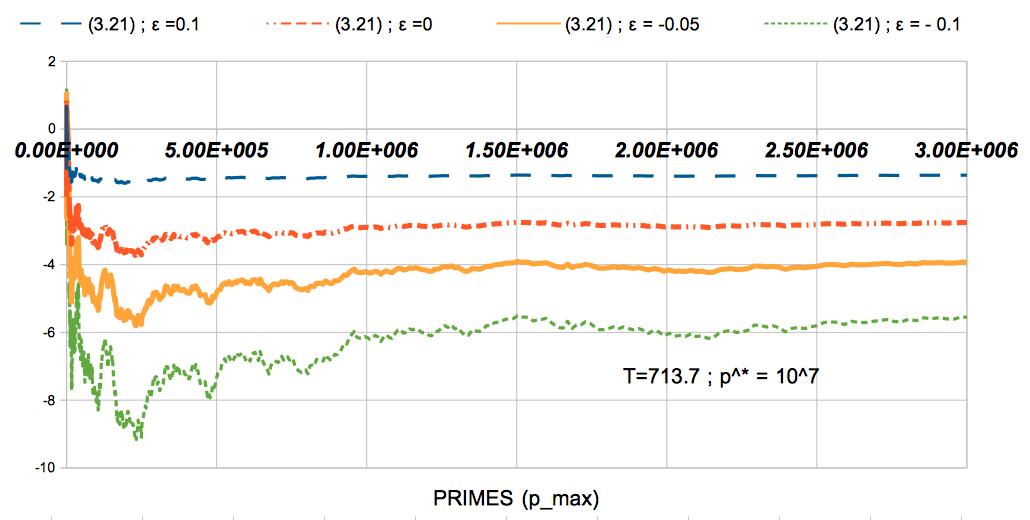} 
\caption{\small {\it   
%ODS  18_06_22ForStabilizationArgument.   sheet  SinSamp 713 7
				Plotted with    (\ref{ProVarPhaseCsi}  ) at fixed T,  and, with  $p_{max} $ from 2   till $3 \times10^6$. We can see  the stabilization trend at increasing primes at different 	 $\epsilon$ values.  % 
				Stabilization is reached at greater primes with lower  algebraic $\epsilon$. Notice $\epsilon $ values from  $-0.1$ to $0.1$%
} }%
\label {StabilizationAtNegativeEps}  
\end{center}
\end{figure}

\begin{figure}[!htbp]
\begin{center}
\includegraphics[width=1.0\textwidth]{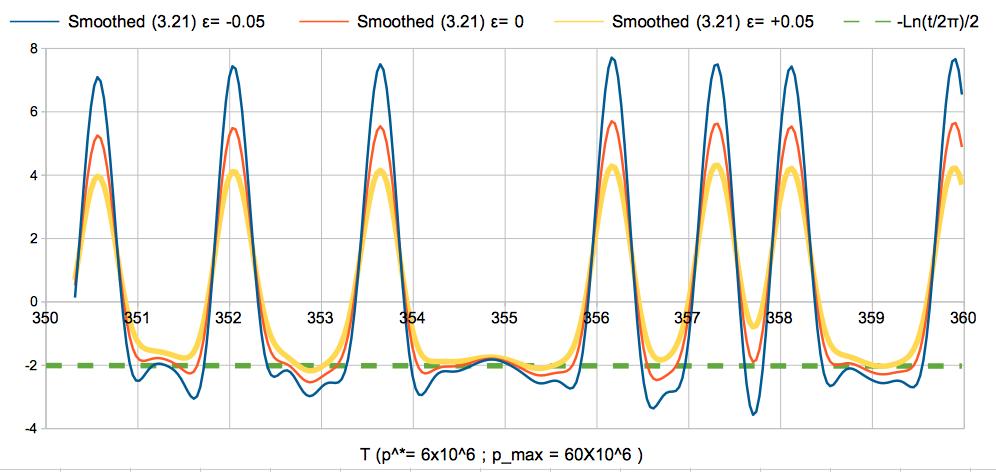} 
\caption{\small {\it   
% ODS  11 01  23  SumMinusInt sheet CorrT350T360   
Computation of  (\ref{ProVarPhaseCsi}  )  with  $p_{max}=60 \times 10^6$  and   $p^*=6 \times 10^6$
 with a sliding window  of length     
 $W= \frac {2 \pi}{\ln(p_{max})}= \frac {2 \pi} {\ln(6\times10^7) } $      
  to counteract Gibbs Phenomenon.
          The shift on t axis is not compensated. 
Notice that, for $|t-t_k^*|>\frac{2 \pi}{\ln(p^*)} $, curves with lower $\epsilon$ are below curve with higher $\epsilon$  as predicted by Lemma 3. If we put  $p_{max}=p^*=6 \times 10^6$ we get a plot   similar to this one, if then  we do smoothing  with $W= \frac {2 \pi}{\ln(p_{max})}= \frac {2 \pi} {\ln(6\times10^6) } $  we get, practically, this same figure. Notice the greater positive peak for negative $\epsilon$. It is due to the discontinuity of $ \Im [ \ln(\zeta(s)(s-1)) ]     $ along $t=t^*_k \  ,  \ \epsilon \le \epsilon^*_k$, here $=0$.
} }
\label {SmthT350T360Pstar6Pmax60}  
\end{center}
\end{figure}

\begin{figure}[]
\begin{center}
\includegraphics[width=1.0\textwidth]{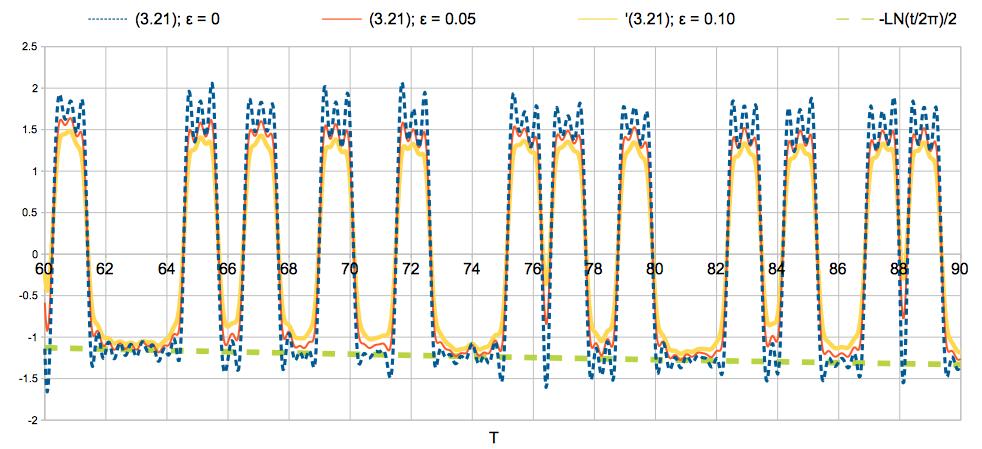}  
\caption{ \small {\it  Here, in interval $60 <t<90$,  it is shown phase variations computed with  (\ref  {ProVarPhaseCsi}) with $p^*=280$  and   $p_{max}=280^3 =21952000$. Used values of $\epsilon$ are : $\epsilon=0 ; +0.05; + 0.1$. %The fat spikes (due to low $p^*$)  produced by $ \zeta(s)$ zeros are compared with  $\frac { \partial \angle [Z(t,\epsilon)] }{\partial t}$ using (\ref {ZSinhECosh1}). 
Far from spikes (i.e. $ \zeta(s)$zeros ) it is apparent the trend  to approach to curve $- \ln\left( \sqrt{\frac{t}{2 \pi}}  \right)$  for $\epsilon=0 $.
%ODS : 02 10 21ZT30T60(+T930+T1830)  T30T90Per3Eps
 % In questo caso e' analizzata la variazione di fase dell'Euler Product della funzione $\zeta(s)$. Qui  $p^*=280$  . Tanto $90/280 \approx 0.32$ e per $ t$  fra 30 e 90 lo spacig fra gli zeri  puo' giustificare un intervallo di media di $2 \pi/ \ln(280) =1.11$. I numeri primi considerati sono quelli fino a $280^3 =21952000$ (cioe' circa 1.3 milioni).  Lontano dagli zeri e' visibile la tendenza ad avvicinarsi alla curva $- \ln\left( \sqrt{\frac{t}{2 \pi}}  \right)$ . I valori di $\epsilon$ usati sono: $\epsilon_0=0 ; \epsilon_1=0.05; \epsilon_2=0.1$
%
%Le oscillazioni ricordano  il Gibbs phenomenon su onda quadra  (\cite  {Fay1999} ).Dove il filtraggio  low-pass nel dominio delle frequenze crea ringing nel dominio del tempo. Qui pero' la situazione  e' piu' complessa.  Non c'e' una sola step-discontinuity' , bensi' molte, e,  di tipo Dirac-delta. %OVVIAMENTE ANCHE QUESTA FIGURA  NON LA MANDEREI ALLA RIVISTA , MA  SU arXiv PERCHE' NO?
}}
\label {T30T90Per3Eps}
\end{center}
\end{figure} 
 
 %   \ref{T60T90xk0}   \ref{T30T90Per3Eps}

\begin{figure}[]
\begin{center}
\includegraphics[width=1.0\textwidth]{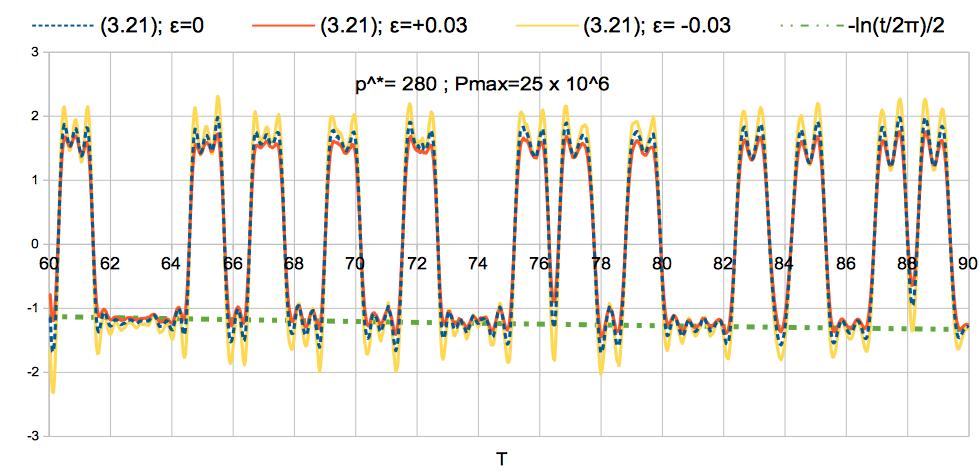}  
\caption{ \small {\it    Same case as fig. \ref {T30T90Per3Eps}  with  $p^*= 280 $, but, with  $p_{max}= 25 \times 10^6$, and with $\epsilon=0 ; + 0.03; -0.03$. The trend toward $- \ln\left( \sqrt{\frac{t}{2 \pi}}  \right)$ far from zeros is apparent for $\epsilon=0 $.
% ODS: 03_09_22  sheet T60T90
%In questo caso $p^*= 280 $ come il fig. \ref{T30T90Per3Eps}  ma qui la sommatoria  (\ref {EulPhaseDerSmoothing2})    e' stata  fatta fino a $x^0_k$ ,cioe'  $p_{max}$ e' dato dalla (\ref {PMaxCheAnnullaIntegr}). Quello indicato in figura e' quello per $t_{min} =60$.
}}
\label {T60T90xk0}
\end{center}
\end{figure}     
 %  {T60T90xk0}  {T30T90Per3Eps}

\clearpage

 %%%%%%%%%%%%%%%%%%%%%%%%%%%%%%%%%%%%%%%%%%%%%%%%%%
\section {  Evaluation of  integral   in(  \ref  {ProVarPhaseCsi})    }
 \label {SumsWithIntegrals} \ \

%            $$ . . .da \ failed \ submit \ 29-01-24  . . . . .$$

 \noindent The aim of this section  is to give a formula for
 the evaluation of  (\ref  {ProVarPhaseCsi}) (reported below for easy reading ):

                                                                                                                                                                                                                                                                                                                                                                                                                                                                                                                                                                                                               \be \label {ProVarPhaseCsi_0}
                                                                                                                                                                                                                                                                                                                                                                                                                                                                                                                                                                                                               \frac{ \Delta \angle[\zeta(s) (s-1)]}{\Delta t} =\frac{ \Delta \Im[\ln(\zeta(s) (s-1))]}{\Delta t} =
                                                                                                                                                                                                                                                                                                                                                                                                                                                                                                                                                                                                               \ee
     
       \be \label {ProVarPhaseCsi_}
 -  \frac {\ln(p^*)} {2 \pi}   \left\{ 
 \left[ \sum_{p =2}^{p_{max} \rightarrow \infty   } \arctan \left( \frac{
 \sin(\ln(p) t      )}{ p^{1/2+\epsilon}-  
\cos( \ln(p) t   )  }\right)\right]_{t-\frac{\pi}{\ln(p^*)}}^{t_2   +\frac{\pi}{\ln(p^*)}   } 
 \right\}  -
\ee

       \be \label {ProVarPhaseCsi1}
  \frac {\ln(p^*)} {2 \pi}   \left\{  - 
 \int_2^  { p_{max} }  
  \frac{2\cos(\ln(x) \ t)}{x^{1/2+ \epsilon}} \sin\left(  \frac{\pi \ln(x)}{\ln(p^*)} \right) 
     d [ Li(x) ] \right\}
\ee  
 
  \noindent In other words we are prompted to compute a definite integral like:

 \be \label  {IntEulPhaseDerSmoothing2TRIS.png}
\frac {\ln(p^*)} {2 \pi}
\int_{2}^{x_k^0} \frac{dp}{\ln(p) \ p^{1/2+\epsilon}} \cos[\ln(p) t] \sin\left( \frac{\pi \ln(p)}{\ln(p^*)}  \right)  = 
  \frac {\ln(p^*)} {2 \pi}  
           \int_2^{ x^0_k} \left\{   \frac{dp} {\ln(p)  } \ \
            \frac{ 
              \sin [\ln(p)t] }{\sqrt{p} p^\epsilon}    \right\}_{t_1=t-\frac{\pi}{\ln(p^*)}}^{t_2=t+ \frac{\pi}{\ln(p^*)}}
  \ee

 % {IntConSinEExp  {IntConCosEExp}
 
\noindent    We can compute it  with  the following indefinite integrals  (\cite [p.~85] {MurrayRSpiegel:2004cz} ):
  
  \be \label {IntConSinEExp}
  \int e^{av} \sin(bv) dv = \frac{e^{av} [a \sin(bv)-b \cos(bv) ] }{a^2 + b^2}
  \ee
  
    \be \label {IntConCosEExp}
  \int e^{av} \cos(bv) dv = \frac{e^{av} [a \cos(bv)+ b \sin(bv) ]}{a^2 + b^2}
  \ee
 
\noindent %  
We  can change integration variable as $ p = e^y$ (i.e  $ y=\ln(p);\quad  dp=e^ydy$), and, reminding that:

 $$
2 \cos(t y) \sin \left(\frac{\pi y}{\ln(p^*)} \right)= \sin\left[ \left(t+\frac{\pi}{\ln(p^*)}  \right)y \right] -
  \sin\left[ \left(t-\frac{\pi}{\ln(p^*)}  \right)y \right] 
 $$
Besides posing : 

\be \label  {iBPiuMenoEA}
 b^\pm= t \pm \frac{\pi}{\ln(p^*)}  \quad \quad and \quad \quad  a= \frac {1}{2}-\epsilon 
 \ee

  \noindent   we can then  write   the integral  in (\ref {ProVarPhaseCsi_} ) as:

 \be \label  {IntEulPhaseDerSmoothing2TRISConY}
  \frac{\ln(p^*)} {2\pi}  
        \int_{y=\ln(s)} ^{y= \ln(p')}
        e^{ay}   [   \sin( b^+ y) -   \sin( b^- y)   ] \frac{dy}{y} 
        \ee

  \noindent  Putting   $u_1(y)=\frac{1}{y}$ and
  
    $$v_1(y)=\frac{e^{ay} [a \sin(b^+y)-b^+ \cos(b^+y) ] }{a^2 +( b^+)^2}-
    \frac{e^{ay} [a \sin(b^-y)-b^- \cos(b^-y) ] }{a^2 +( b^-)^2}
    $$

    with integration by parts we get:

 $$
  \frac{\ln(p^*)} {2\pi}  
        \int_{y=\ln(s)=y_1} ^{y= \ln(p')=y_2}
        e^{ay}   [   \sin( b^+ y) -   \sin( b^- y)   ] \frac{dy}{y}  =   [  u_1v_1]_{y_1}^{y_2} -  \int_{y_1}^{y_2}v_1 du_1  =
 $$  
  
   \be \label {PrimIter}
   =\left\{ \frac{e^{ay} }{y} \left[ \frac{[a \sin(b^+y)-b^+ \cos(b^+y) ] }{a^2 +( b^+)^2}-
    \frac{ [a \sin(b^-y)-b^- \cos(b^-y) ] }{a^2 +( b^-)^2} \right]   \right\}_{y_1}^{y_2}+
    \ee
    $$
    \int_{y_1}^{y_2} \frac{e^{ay} }{y^2}\left[ \frac{[a \sin(b^+y)-b^+ \cos(b^+y) ] }{a^2 +( b^+)^2}-
    \frac{ [a \sin(b^-y)-b^- \cos(b^-y) ] }{a^2 +( b^-)^2} \right] dy
   $$

  \noindent   Putting now 
  
        $$A^\pm_1=A^\pm=\frac{a}{a^2+(b^\pm)^2}  \quad \quad and \quad \quad  B^\pm_1= B^\pm=\frac{b^\pm}{a^2+(b^\pm)^2} $$

  \noindent   we can write
     
      $$     u_1(y)=\frac{1}{y}       \quad \quad and \quad \quad  v_1(y) =
         e^{ay}  \{  A_1^+\sin(b^+y) - B_1^+\cos(b^+y)
          -A_1^-\sin(b^-y) + B_1^-\cos(b^-y) \}                                                          $$

     \noindent  so  (\ref   {PrimIter}) becames:
     
      \be \label {PrimIter2}
       [  u_1v_1]_{y_1}^{y_2}+ \int_{y_1}^{y_2} \frac{e^{ay} }{y^2}\left[ 
      A_1^+\sin(b^+y) - B_1^+\cos(b^+y)
          -A_1^-\sin(b^-y) + B_1^-\cos(b^-y)
    \right] dy
      \ee

  The integral in  (\ref   {PrimIter2}) can be computed by iteration  using  (\ref {IntConSinEExp}) and
     (\ref  {IntConCosEExp}) 
     
   \noindent    and for $j >1$:
     
     $$  A_j^\pm=A^\pm A^\pm_{j-1}  -   B^\pm B^\pm_{j-1}     \quad \quad and \quad \quad
        B_j^\pm=A^\pm B^\pm_{j-1}  -   B^\pm A^\pm_{j-1}                 $$
        
    \noindent      with:

     \noindent       for $j>1$

          $$ u_j=-\frac {\partial   u_{j-1}(y)    } {\partial y}   \quad \quad  ; \quad \quad   \frac {\partial v_{j}  (y) } {\partial y}  =v_{j-1}(y)$$

            \be \label {Iteration}     u_j(y)=\frac{(j-1)!}{y^j}       \quad \quad and \quad \quad 
             v_{j}(y) =
        e^{ay}  \{  A_j^+\sin(b^+y) - B_j^+\cos(b^+y)
          -A_j^-\sin(b^-y) + B_j^-\cos(b^-y) \}                                                          \ee

 Of course the presence of factorial forces to block  $J_{max}$ to a suitable value (so it is an approximate solution). A practical bond is  $ J << \lfloor T \rfloor$ as { \bf  $A_J^\pm$ and $B_J^\pm$  decrease like  $T^{-J}$} .%

 $$ \frac{\ln(p^*)} {2\pi}  
        \int_{y=\ln(s)=y_1} ^{y= \ln(p')=y_2}
        e^{ay}   [   \sin( b^+ y) -   \sin( b^- y)   ] \frac{dy}{y}  =   [  u_1v_1]_{y_1}^{y_2} -  \int_{y_1}^{y_2}v_1 du_1  =$$
        
   \be \label {SolApprox} \sim  \sum_{j=1}^{j=J_{max}}   [  u_j v_j]_{y=\ln(s)}^{y=\ln(p')} \quad J_{max} << t  \ee

\noindent By known integrals  $\int e^{av} \sin(bv) dv$ and  $\int e^{av} \cos(bv) dv$, we have :
$$ \frac{\ln(p^*)} {2\pi}  
        \int_{y_1} ^{y_2}
        e^{ay}   [   \sin( b^+ y) -   \sin( b^- y)   ] \frac{dy}{y} =
        $$

        \be \label {IntByPartsFirst}
        \frac{\ln(p^*)}{2 \pi} \left\{ 
        \frac{e^{ay}} {y}  \left[ 
         \frac{   [a \sin(b^+y)-b^+ \cos(b^+y) ] }{a^2 +( b^+)^2}-
         \frac{
     [a \sin(b^-y)-b^- \cos(b^-y) ] }{a^2 +( b^-)^2} 
      \right] 
      \right \}_{y_1}^{y_2} + \ee

      $$\frac{\ln(p^*)}{2 \pi} 
      \int_{y_1}^{y_2}
        \frac{e^{ay }\ dy} {y^2}  \left[ 
         \frac{   [a \sin(b^+y)-b^+ \cos(b^+y) ] }{a^2 +( b^+)^2}-
         \frac{
     [a \sin(b^-y)-b^- \cos(b^-y) ] }{a^2 +( b^-)^2} 
      \right]$$

\

In  (\ref {IntByPartsFirst}) we can choose for $y_1$ whatever point. In particular  we choose   a point  $y_1$ with $ \left[ 
         \frac{   [a \sin(b^+y_1)-b^+ \cos(b^+y_1) ] }{a^2 +( b^+)^2}-
         \frac{
     [a \sin(b^-y_1)-b^- \cos(b^-y_1) ] }{a^2 +( b^-)^2} 
      \right]=0$, so  (\ref {IntByPartsFirst})  becomes:

\be \label {DeveEssereZero}   \frac{\ln(p^*)}{2 \pi}  \left\{ 
        \frac{e^{ay_2}} {y_2}\right\}
 \left[ 
         \frac{   [a \sin(b^+y_2)-b^+ \cos(b^+y_2) ] }{a^2 +( b^+)^2}-
         \frac{
     [a \sin(b^-y_2)-b^- \cos(b^-y_2) ] }{a^2 +( b^-)^2} 
      \right] + \ee
      
    $$ \frac{\ln(p^*)}{%
    2 \pi}
   \int_{y_1}^{y_2}
        \frac{e^{ay }\ dy} {y^2}  \left[ 
         \frac{   [a \sin(b^+y)-b^+ \cos(b^+y) ] }{a^2 +( b^+)^2}-
         \frac{
     [a \sin(b^-y)-b^- \cos(b^-y) ] }{a^2 +( b^-)^2} 
      \right]
      $$
      
  \noindent    If we drop the integral in  (\ref {DeveEssereZero}  ) we have (\ref {SolApprox}) for $J_{max} =1$.  %In previous and follo computations 
     Everywhere we have used $J_{max}=3$.

% $$ . . . . .$$

\subsection{Numerical example on convergence of (\ref  {ProVarPhaseCsi}  )   }

\label  {NumEx}

\

% \subsection {Computation of phase variation through Euler product} \label {EulPrPhaVar}

\noindent Convergence of  (\ref  {ProVarPhaseCsi}  ),  or   (\ref{ProVarPhaseCsi_0})  (the same), stems from convergence of  $\zeta(s)(s-1)$, see subsection \ref {CONJ}, but, it is worthwhile  to verify how it happens numerically.

% Nelle figure  sono calcolati i pezzi della (\ref {SpetrtroFlutt_Rho1mezz}) con i campionamenti  dati dalla  (\ref{TransizMenoPiMezzi__} )  e (\ref{TransizPiuPiMezzi__} ) .

\noindent Expression  (\ref{ProVarPhaseCsi_0})  i.e.  (\ref{ProVarPhaseCsi})    is expressed like a sum over primes till $p_{max}$, i.e.  (\ref  {ProVarPhaseCsi_}), minus an integral  on $d[Li(x)]$, (x real  prime) always till $p_{max}$ i.e.  (\ref  {ProVarPhaseCsi1}). Here we want to see how this difference behaves computationally. The funtions under the sum and in the integral are respectively.

\begin{itemize}

\item

  \noindent For the sum  a function using $y$,  as continuous variable, in place of discrete $p$:  
    \be \label {ThoughtContinuous} 
    F(y,\epsilon)= \left\{ \arctan \left( \frac{
 \sin(\ln(y) t     )}{ p^{1/2+\epsilon}-  
\cos( \ln(y) t   )  }\right) \right \}_{t-\frac{\pi}{\ln(p^*)}   }^{t+\frac{\pi}{\ln(p^*)}  }
\ee   
      
 \noindent   $ y_{0t}(k+1), y'_{0t}(k)$ are defined,for $x < p^*$,  as positive slope zero transitions   and negative slope zero transition respectively  (similarly to    $ x_{0t}(k+1), x'_{0t}(k)$ in  (\ref{TransizMenoPiMezzi__} ) and (\ref{TransizPiuPiMezzi__} )). See    \cite [p.~10]  {Giovanni Lodone Nov2024} .

 \item
 \noindent For the integral the function is (see \cite  {Giovanni Lodone Nov2024} appendix B fig 6) :

 \be \label {SingleOillLi1}  
 %
     %
 % \frac {\ln(p^*)} {\pi \phi(q)}  
       %    \int_{ x >x_{0t} (k,h)  \ : \  p \equiv h (mod \  q)}^{ x <x'_{0t} (k,h)} 
            \frac{ 
            %\ 
              \cos [\ln(x)t] ]\sin [\pi \ln(x) /\ln(p^*)] }{\sqrt{x} x^\epsilon}   dLi(x)
       %       \right|
 \ee

\end{itemize}

\noindent  For $p$ or $y$  big, we have : $ x_{0t}(k+1) \rightarrow  y_{0t}(k+1)$ and $ x'_{0t}(k) \rightarrow  y'_{0t}(k)$, and, besides, we have a relative maximum of the sum or the integral  along $x$ or $y$ at $ x'_{0t}(k) \rightarrow  y'_{0t}(k)$, (\ref{TransizPiuPiMezzi__} ),  because the sum or the integral  stops increasing and begins decreasing. While we have a relative minimum of the sum or the integral  along $x$ or $y$ at $ x_{0t}(k) \rightarrow  y_{0t}(k)$,  (\ref{TransizMenoPiMezzi__} ),  because the sum or the integral  stops decreasing and begins increasing.
In following figures are computed in order  :
\begin{itemize}
\item  the envelop of maxima of  (\ref{ProVarPhaseCsi_})  sampling at     (\ref{TransizPiuPiMezzi__} )       and the envelope of minima sampling at   (\ref{TransizMenoPiMezzi__} ), for $(p^*)^m<x<(p^*)^{m+1}$ with even $m$, ($m=0$ in figure). The reverse holds  with odd $m$ (not shown). The varying gap between maxima and minima envelopes increases with $m$, but it is sqeezed (close to zero) at $p=(p^*)^{m+1}$.

\item  the envelop of maxima of  (\ref{ProVarPhaseCsi1})  sampling at     (\ref{TransizPiuPiMezzi__} )       and the envelope of minima sampling at   (\ref{TransizMenoPiMezzi__} ), for $(p^*)^m<x<(p^*)^{m+1}$ with even $m$, ( $0$ in figure). The reverse holds  with odd $m$ (not shown). The varying gap between maxima and minima envelopes increases with $m$, but it is sqeezed (sometimes to  zero) at $p=(p^*)^{m+1}$.

\item  the difference  :   (\ref{ProVarPhaseCsi_})  minus  (\ref{ProVarPhaseCsi1}), sampling at   (\ref{TransizMenoPiMezzi__} ) and  at  (\ref{TransizPiuPiMezzi__} )
\end{itemize}

\begin{figure}[]
\begin{center}
\includegraphics[width=1.0\textwidth]{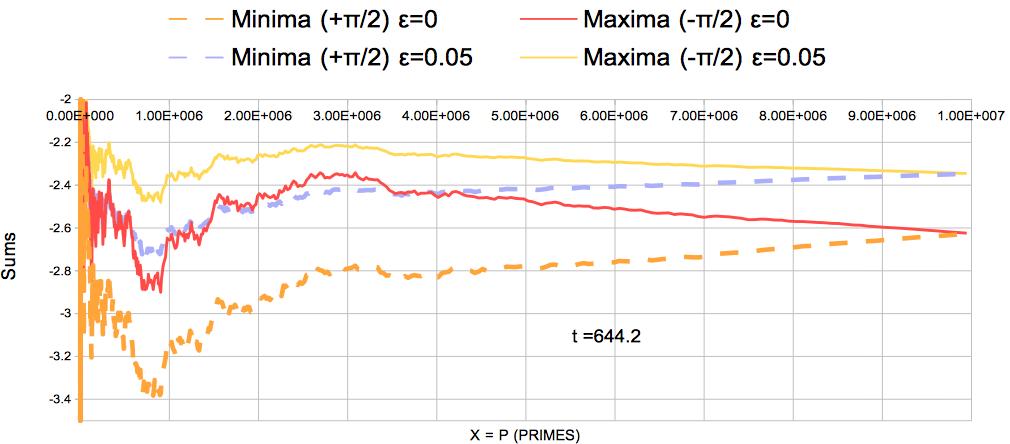}  
\caption{ \small {\it   
% ODS: 23_01_22
Envelope of minima for sum 
%Sum in    (\ref {ProVarPhaseCsi}) 
in expression    (\ref {ProVarPhaseCsi_})
%Sums  taken alone in (\ref {ProVarPhaseCsi}) , i.e $  -  \frac {\ln(p^*)} {2 \pi}   \left\{ 
% \left[ \sum_{p =2}^{p  } \arctan \left( \frac{ \sin(\ln(p) t      )}{ p^{1/2+\epsilon}-  \cos( \ln(p) t   )  }\right)\right]_{t-\frac{\pi}{\ln(p^*)}}^{t   +\frac{\pi}{\ln(p^*)}   }  \right\}  $
, sampled at                        increasing  cosine  $y_{0t}<p^*$  (\ref{TransizMenoPiMezzi__} ). And envelope of 
, local maxima of oscillation,  at decreasing  cosine $y'_{0t} <p^*$   (\ref {TransizPiuPiMezzi__}). % local maximum of oscillations.  To compare them the abscissa for maxima and minima in same oscillation is the same:  $\frac {x'_{0t} +x_{0t}}{2}$. 
We choose in expression   (\ref {ProVarPhaseCsi_}) $p^*\approx 10^7$,
% (i.e.   ( \ref {EulPhaseDerSmoothing2}) ) 
 for $t=644.2$ and $\epsilon = 0 $ or $\epsilon= 0.05$. Dashed lines (relative minima) are with positive   slope cosine at zero crossing  sampling (\ref {TransizPiuPiMezzi__} ) $ + \pi/2$ in exponent. Continuous  lines ( relative maxima) are with  negative         slope cosine at zero crossing  sampling (\ref {TransizMenoPiMezzi__} ) $ - \pi/2$ in exponent. Notice that  (\ref {TransizPiuPiMezzi__} )  and  (\ref {TransizMenoPiMezzi__} )  are defined on $\int dLi(x)$ integral, but for $p$, i.e. $y$, big they coincide: i.e.  $y'_{0t} \rightarrow x'_{0t}$ and $y_{0t} \rightarrow x_{0t}$.
}}
\label {SumsT644}
\end{center}
\end{figure}     
 
  %     \ref{T140T160PstarSmaller}   . . . \ref{IntegralsT644}

\begin{figure}[]
\begin{center}
\includegraphics[width=1.0\textwidth]{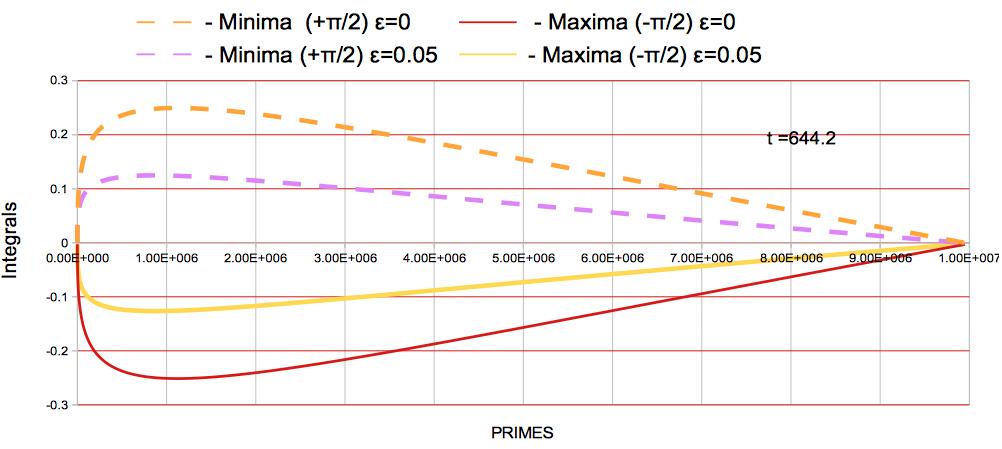}  
\caption{ \small {\it   
% ODS: 23_01_22
Envelopes of  minus relative minima and minus  relative maxima of   expression(\ref {ProVarPhaseCsi1}) 
  for $t=644.2$ and $\epsilon = 0 $ or $\epsilon= 0.05$. 
     Dashed lines are  with positive  slope zero crossing  cosine  sampling (\ref {TransizPiuPiMezzi__} ) . Continuous  lines are with negative  slope zero crossing     cosine  sampling (\ref {TransizMenoPiMezzi__} ).  Note the minus sign. For $(p^*)^m<x<(p^*)^{m+1}$ the oscillations increase in amplitude with sign change with $m$ changing parity.
}}
\label {IntegralsT644}
\end{center}
\end{figure}     
 
%  {SpetrtroFlutt2_}

\begin{figure}[]
\begin{center}
\includegraphics[width=1.0\textwidth]{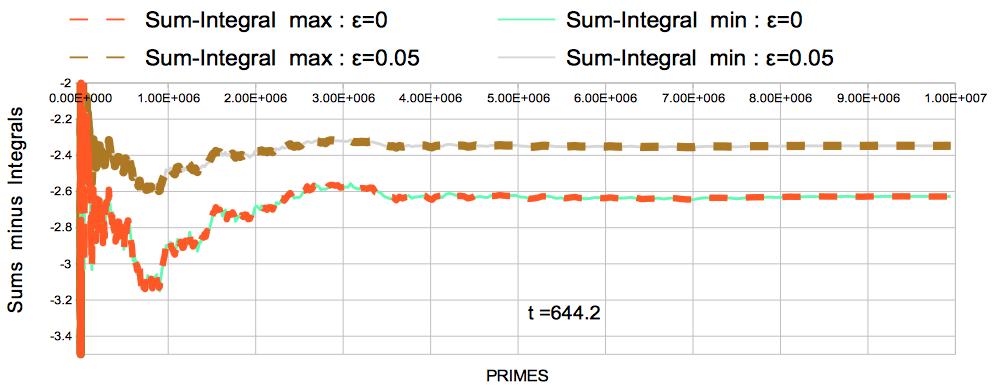}  
\caption{ \small {\it   
% ODS: 23_01_22
Sum   (\ref {ProVarPhaseCsi_}) $ -  Integral$  (\ref {ProVarPhaseCsi1}), sampled at $x_{0t}$ and $x'_{0t}$.See  (\ref{TransizPiuPiMezzi__} )       and (\ref{TransizMenoPiMezzi__} ).
%SumsMinusIntegrals $\left[   \frac {\Delta \angle [ \zeta_{EP}(s)] }{\Delta t}     \right]_{\sum_{\forall p<p'<p^*}} - \left[    \frac {\partial   \angle[ EulerProduct]  }{\partial  t}      \right]_{\int_2^{p'<p^*} dp} $ (i.e.   ( \ref {EulPhaseDerSmoothing2})  + (\ref {IntEulPhaseDerSmoothing2TRISConY}) = ( \ref {SpetrtroFlutt2_}) )  for $t=644.2$ and $\epsilon = 0 $ or $\epsilon= 0.05$. Dashed lines are with  negative slope   at   zero crossing cosine  sampling (\ref {TransizPiuPiMezzi__} ) . Continuous  lines are with positive slope cosine at zero crossing  sampling (\ref {TransizMenoPiMezzi__} ). It is apparent that  ( \ref {SpetrtroFlutt2_}) is invariant with respect  to the chosen % slope of  zero crossingsampling. {\bf Note that the last sampling point, close to $p^*$, is congruent with previous samplings}.%Note that the constant level reached at high $p$  is midway from curves  in fig, \ref  {SumsT644}.
Maxima and minima envelops shrink in a single curve. Notice  that as we are far from zeros  the limiting value at $\epsilon=0.05$ is greater  than limiting value at $\epsilon=0$
}}
\label {SumsMinusIntegralsT644}
\end{center}
\end{figure}

   \noindent  It is apparent that  (\ref {ProVarPhaseCsi_}), and  (\ref {ProVarPhaseCsi1}) oscillates between whatever large positive and negative values.  But  (\ref {ProVarPhaseCsi_}) $ - $  (\ref {ProVarPhaseCsi1}), after an initial transition lasting ( it depends on $t$ and $\epsilon$), in this case,  till primes less then $3$ or $4$ $\times 10^6$, %sum  (\ref {ProVarPhaseCsi_}) minus integral (\ref {ProVarPhaseCsi1}) 
   it stabilizes  afterward  to a converging value. The same value both with maxima difference and minima difference.

\section{
Approximation to   (\ref{EulPhaseDerSmoothing2})
}
The difference between  the infinite sum (\ref{EulPhaseDerSmoothing2})  and  the % integral in (\ref{Trasf1ImPartDerivSuT}) , i.e. (\ref {ProVarPhaseCsi}),  
approximation : 
\be \label {Approx}
 \frac {\ln(p^*)}{2 \pi} 
 \sum_{  p<p_{max}      } \frac{2 \cos [\ln(p) t )] ] \sin\left( \frac{\ln(p) \pi}{\ln(p^*)}   \right) }{p^{1/2+\epsilon}}
\ee
\noindent obtained by:

\be \label {FdipparamChiH}
 \left\{ \arctan \left( \frac{
 \sin(\ln(p) t   
  )}{ p^{1/2+\epsilon}-  
\cos( \ln(p) t )  }\right) \right \}_{t-\frac{\pi}{\ln(p^*)}   }^{t+\frac{\pi}{\ln(p^*)}  } \rightarrow
\ee
$$\left\{  \left( \frac{
 \sin(\ln(p) t     )}{ p^{1/2+\epsilon}-  
\cos( \ln(p) t  )  }\right) \right \}_{t-\frac{\pi}{\ln(p^*)}   }^{t+\frac{\pi}{\ln(p^*)}  } \rightarrow
\left\{  \left( \frac{
 \sin(\ln(p) t    )}{ p^{1/2+\epsilon}   }\right) \right \}_{t-\frac{\pi}{\ln(p^*)}   }^{t+\frac{\pi}{\ln(p^*)}  } =$$

$$=  \frac{ \sin \left(  \ln(p) \left\{ t +\frac{\pi}{\ln(p^*)} \right\}      \right) -
  \sin \left( \ln(p) \left\{ t - \frac{\pi}{\ln(p^*)} \right\}  \right)}
 {p^{1/2+\epsilon}}  =
 $$

% (ref {DerLAngleSuTL})   (ref {FdaLgaritmicInt} )

%tende a 
\be \label {FdaLgaritmicInt}
 \frac{2 \ \cos( \ln(p) t   )   \sin\left( \frac{\pi \ln(p)}{\ln(p^*)} \right)  }  
  {p^{1/2+\epsilon}}  
  \ee

is convergent for $\epsilon > 0$.

\noindent PROOF. Using Taylor formula for arctangent the difference between  (\ref{EulPhaseDerSmoothing2}))  and (\ref {Approx}),  using (\ref  {FdipparamChiH}) or (\ref {FdaLgaritmicInt}),   can be written as:

$$
-  \left. \frac {\ln(p^*)} {2 \pi}   \sum_{ p<p_{max}      } \arctan \left( \frac{
 \sin(\ln(p) t     )}{ p^{1/2+\epsilon}-  
\cos( \ln(p) t   )  }\right)\right|_{t_1}^{t_2}
+$$
$$
 \frac {\ln(p^*)}{2 \pi} 
 \sum_{  p<p_{max}    } \frac{2 \cos [\ln(p) t  ] \sin\left( \frac{\ln(p) \pi}{\ln(p^*)}   \right) }{p^{1/2+\epsilon}} =
$$

\be \label {SumInPandN}
 -  \left. \frac {\ln(p^*)} {2 \pi}  \ \ \  \sum_{ p<p_{max}     } \  \ \ \ \  \sum_{n  \ odd >1} \frac{(-1)^{(n-1)/2}}{n}
  \left[   \left( \frac{
 \sin(\ln(p) t     )}{ p^{1/2+\epsilon  } %-  
- \cos(\ln(p) t     )  }\right)^n   \right] \right |_{t_1}^{t_2}
- 
\ee

\be \label {SumPminus2}
\left. \frac {\ln(p^*)} {2 \pi}   \sum_{ p<p_{max}      }  \left( \frac{
 \sin(\ln(p) t     )    \cos(\ln(p) t     )   }{ [p^{1/2+\epsilon}-  
\cos( \ln(p) t  ) ]    p^{1/2+\epsilon}   }\right)\right|_{t_1}^{t_2}
\ee

 \noindent Developing last sum we have, focusing on   main term (i.e. (\ref{SumPminus2}):

\be \label {SumPminus3}
\left .\frac {\ln(p^*)} {4 \pi}   \sum_{ p<p_{max}     }  \left( \frac{
 \sin( 2( \ln(p) t     ) )     }{ [p^{1/2+\epsilon}-  
\cos( \ln(p) t  ) ]    p^{1/2+\epsilon}   }\right)\right|_{t_1}^{t_2} \approx
\ee

\be \label {SumPminus4}
 \frac {\ln(p^*)} {4 \pi}   
\sum_{ p<p_{max}      } 
 \left( \frac{
 \sin( 2( \ln(p) t_2    ) )  -
 \sin( 2( \ln(p) t_1     ) ) 
    }
 { [p^{1/2+\epsilon}-  
\cos( \ln(p) t_2   ) ]    [p^{1/2+\epsilon}-  
\cos( \ln(p) t_1   ) ]     }\right)  
\ee

\noindent It is apparent that (\ref{SumPminus4}), for $p_{max}\rightarrow \infty$, it  is absolutely convergent  from $\epsilon>0$ . So for $\epsilon>0$ the difference (\ref  {SumInPandN}) is always limited.

 END OF PROOF

%    \end{comment}

\section {Expression (\ref{ProVarPhaseCsi})  % irrespective of the meanning that could connect it, 
 %to $\left\{  \frac{ \Delta \Im [\ln \{ \zeta(s)(s-1) \}]}{\Delta t} \right\}_{p^*, p_{max}} $ for $\epsilon>1/2$,
  is convergent for  $\epsilon >0$ }

%%%%%%%%%%%%%%%%%%%%%%%%%%%%%%
 
% \subsubsection{
% Lemma 2.A :
% { \bf Expression (\ref{ProVarPhaseCsi}), irrespective of the meanning that could connect it, 
 %to $\left\{  \frac{ \Delta \Im [\ln \{ \zeta(s)(s-1) \}]}{\Delta t} \right\}_{p^*, p_{max}} $ for $\epsilon>1/2$,
 % is convergent for  $\epsilon >0$}.
% }

 \ 
 
\noindent For easy reading we write again expression (\ref{ProVarPhaseCsi}) 
\label {Conv3_21}
 
 $$ -  \frac {\ln(p^*)} {2 \pi}   \left\{ 
 \left[ \sum_{p =2}^{p_{max}   } \arctan \left( \frac{
 \sin(\ln(p) t      )}{ p^{1/2+\epsilon}-  
\cos( \ln(p) t   )  }\right)\right]_{t-\frac{\pi}{\ln(p^*)}}^{t   +\frac{\pi}{\ln(p^*)}   } - 
 \int_2^  { p_{max} }  
  \frac{2\cos(\ln(x) \ t)}{x^{1/2+ \epsilon}} \sin\left(  \frac{\pi \ln(x)}{\ln(p^*)} \right) 
     d [ Li(x) ] \right\}$$
 
 Using (\ref {Approx}), in oder to look  only for divergencies, we can write:
 \be \label {ApproxForDiv}
  -  \frac {\ln(p^*)} {2 \pi}   \left\{ 
 \left[ 
 \sum_{p =2}^{p_{max}   }
 \frac{2 \cos [\ln(p) t )] ] \sin\left( \frac{\ln(p) \pi}{\ln(p^*)}   \right) }{p^{1/2+\epsilon}}
 % \arctan \left( \frac{ \sin(\ln(p) t      )}{ p^{1/2+\epsilon}-  
%\cos( \ln(p) t   )  }\right)
\right]  %_{t-\frac{\pi}{\ln(p^*)}}^{t   +\frac{\pi}{\ln(p^*)}   }
 - 
 \int_2^  { p_{max} }  
  \frac{2\cos(\ln(x) \ t)}{x^{1/2+ \epsilon}} \sin\left(  \frac{\pi \ln(x)}{\ln(p^*)} \right) 
     d [ Li(x) ] \right\}
 \ee
 
 \noindent we will show that (\ref {ApproxForDiv}) converges $\forall t \ , \ \forall \epsilon>0$
 
  \noindent PROOF  In each interval containing only the prime p  the discrete sum  can be written as a defined integral, and, the infinite integral in $y$ can be split in an infinite sum of defined integrals.
 \noindent Supposing the prime $p$ (fixed in each interval, while $t$ is fixed everywhere)  is surrounded by neighbours primes at $p- 2\alpha \ln(p) < p <p + 2\beta \ln(p) $ 
 
 %\footnote{ If the  distance among primes was exactly the mean distance, then $\alpha=\beta=1/2$. This is not the case so $\alpha$  and $\beta$ manage the practical cases. Notice a potential application to $h-$class primes $mod \ q$ ( see \cite  {Giovanni Lodone Nov2024}). The interval with mean $\alpha,\beta = 1/2$ would become:  $$ p_{h-class} -\phi(q) \alpha \ln( p_{h-class} ) < y <p_{h-class} +\phi(q) \beta \ln( p_{h-class} )$$  This  could be perhaps used to associate phase variation peaks to particular $h-$classes of primes }
 ,  we can write  using approximation  (\ref{Approx}):

 \be \label {ApprMinus0_1}
 -\frac {\ln(p^*)}{ \pi} 
 \sum_{  p<p_{max}     } 
 \int_{p- \alpha \ln(p) }^{p + \beta \ln(p)}  
 \ee
 $$
 \left\{  
 \frac{ \cos [\ln(p) t  ] \sin\left( \frac{\ln(p) \pi}{\ln(p^*)}   \right) } { (\alpha+\beta) \ln(p) p^{1/2+\epsilon}}
  -
   \frac{\cos(\ln(y) \ t    )} {y^{1/2+ \epsilon}} \sin\left(  \frac{\pi \ln(y)}{\ln(p^*)} \right)  \frac {1}{\ln(y)}  \right\} dy =
  $$

   \be \label {ApprMinus0_2}
 -\frac {\ln(p^*)}{ \pi} 
 \sum_{  p<p_{max}     } 
 \int_{p- \alpha \ln(p) }^{p + \beta \ln(p)}  
 \ee
 $$
 \left\{  
 \frac{
  y^{1/2+ \epsilon}  \ln(y)      \cos [\ln(p) t ] \sin\left( \frac{\ln(p) \pi}{\ln(p^*)}   \right)      -  (\alpha+\beta) \ln(p) p^{1/2+\epsilon}  \cos(\ln(y) \ t )   \sin\left(  \frac{\pi \ln(y)}{\ln(p^*)} \right)  
  } 
 { (\alpha+\beta) \ln(p) p^{1/2+\epsilon}    y^{1/2+ \epsilon}  \ln(y) }
   \right\} dy 
  $$
  \noindent Inside the interval  $p- \alpha \ln(p) < y <p + \beta \ln(p) $ ( $ \Delta y= y-p$) ,  
  to first order we have:
  $$ \cos(\ln(y) \ t )  =    \cos(\ln(p) \ t ) -\frac{t(y-p)}{p}\sin( \ln(p) \ t )$$
 
  $$ \sin\left(  \frac{\pi \ln(y)}{\ln(p^*)} \right) =
    \sin\left(  \frac{\pi \ln(p)}{\ln(p^*)} \right) + \frac{\pi (y-p)}{ p \ln(p^*)}    \cos\left(  \frac{\pi \ln(p)}{\ln(p^*)} \right)$$
 
 so { \bf at first order } in  $\frac{(y-p)}{p}$, with $p >> Max(\alpha, \beta) \times \ln(p)$:
 
 $$  y^{1/2+ \epsilon}  \ln(y)      \cos [\ln(p) t  ] \sin\left( \frac{\ln(p) \pi}{\ln(p^*)}   \right)      -  (\alpha+\beta) \ln(p) p^{1/2+\epsilon}  \cos(\ln(y) \ t  )   \sin\left(  \frac{\pi \ln(y)}{\ln(p^*)} \right)  =
 $$
 $$
 \cos(\ln(p) \ t ) \sin\left(  \frac{\pi \ln(p)}{\ln(p^*)} \right)
  [ y^{1/2+ \epsilon}  \ln(y) -   (\alpha+\beta) \ln(p) p^{1/2+\epsilon} ]  +
 $$
 $$
 \frac{(y-p)}{p} \left\{   \frac{\pi }{  \ln(p^*)}    \cos\left(  \frac{\pi \ln(p)}{\ln(p^*)} \right) 
 -t\sin( \ln(p) \ t  )
  \right\}
 $$

 So the numerator of the fraction under integral in (\ref{ApprMinus0_2}) can be split in two parts resulting in two integrals:

  \be \label {FirstO_} 
  \frac {\ln(p^*)}{ \pi}  
  \sum_{  p<p_{max}   \ : \ gcd(p,q)=1   } 
 \int_{p- \alpha \ln(p) }^{p + \beta \ln(p)}  
  \frac{   \frac{\pi }{  \ln(p^*)}    \cos\left(  \frac{\pi \ln(p)}{\ln(p^*)} \right) 
 -t\sin( \ln(p) \ t   ) } 
  { (\alpha+\beta) \ln(p) p^{1/2+\epsilon}    y^{1/2+ \epsilon}     \ln(y) } \ 
    \frac{(y-p)}{p} dy   
    \ee
 
 that in absolute value is less than
 
  $$
  \frac {\ln(p^*)}{ \pi}   \sum_{  p<p_{max}      } 
 \int_{p- \alpha \ln(p) }^{p + \beta \ln(p)}  
  \frac{   t +1} 
  { (\alpha+\beta) \ln(p) p^{1/2+\epsilon}    y^{1/2+ \epsilon}     \ln(y) } \ 
    \frac{(\alpha+\beta) \ln(p) }{p} dy   = $$
    $$
  \frac {\ln(p^*)}{ \pi}     \sum_{  p<p_{max}      } 
 \int_{p- \alpha \ln(p) }^{p + \beta \ln(p)}  
  \frac{   t+1 } 
  {  p^{3/2+\epsilon}    y^{1/2+ \epsilon}     \ln(y) } \ 
    dy   
     $$
 
 Wich is surely  convergent for $p_{max}\rightarrow \infty$  ,  $\forall t$  for $\epsilon>0$  (and beyond).
 
 And

  $$\frac {\ln(p^*)}{ \pi}   \sum_{  p<p_{max}     } 
 \int_{p- \alpha \ln(p) }^{p + \beta \ln(p)}  
  \frac{  \cos(\ln(p) \ t   ) \sin\left(  \frac{\pi \ln(p)}{\ln(p^*)} \right)
  [ y^{1/2+ \epsilon}  \ln(y) -   (\alpha+\beta) \ln(p) p^{1/2+\epsilon} ]  } 
  { (\alpha+\beta) \ln(p) p^{1/2+\epsilon}    y^{1/2+ \epsilon}  \ln(y) }
   dy   $$
 
\noindent Now 
 $$\left|    \cos(\ln(p) \ t ) \ \frac {\ln(p^*)}{ \pi}  \  \sin\left(  \frac{\pi \ln(p)}{\ln(p^*)}\right) \right | \le 1$$ is upper bouded by 1,
 and 
 
 $y^{1/2+ \epsilon}  \ln(y) -   (\alpha+\beta) \ln(p) p^{1/2+\epsilon} $ in the interval  is upper bounded by 
 $(\alpha+\beta) \ln(p)    \ln(y) $. 
 
 Because,  in interval  $p- \alpha \ln(p) < y <p + \beta \ln(p) $, 
 
 we have:
 
   $|y^{1/2+ \epsilon}   -    p^{1/2+\epsilon} | <\ln(y) < Max(,\alpha,\beta) \ln(p) \ \ ; \forall \epsilon:  0<\epsilon \le 1/2 $ 
 
 but 
 
  $ |  \ln(y) -   (\alpha+\beta) \ln(p) |   <  (\alpha+\beta) \ln(p) $

  so 
  
  $y^{1/2+ \epsilon}  \ln(y) -   (\alpha+\beta) \ln(p) p^{1/2+\epsilon}  <(\alpha+\beta) \ln(p)    \ln(y) $ 
  
%  DA RIVEDERE !!!!!!

\noindent Orders greater than the first    (i.e. all terms but the first in (\ref {LogEP2}) )   lead to  convergent integrals  like their sum ( at least for $\epsilon>0$) like (\ref{FirstO_}).

 \noindent  So (\ref  {ApprMinus0_1}) in absolute value  is upper  bounded by 
 
  \be \label {ApprMinus0_3}
 % \frac {\ln(p^*)}{ \pi} 
 \sum_{  p<p_{max}     } 
 \int_{p- \alpha \ln(p) }^{p + \beta \ln(p)}  
 \left\{  
 \frac{1 } 
 {  p^{1/2+\epsilon}    y^{1/2+ \epsilon}   }
   \right\} dy 
  \ee

  \noindent This means that for $\epsilon>0$ it is absolutely convergent.
  
  \noindent So , by an argument identical to   Lemma2.A,   also 
  (\ref {ProVarPhaseCsi}), is absolutely  convergent, for $p_{max} =(p^*)^m \rightarrow \infty \, m=1,2, 3 . .$ \footnote {This relation between $ p^*$ and $p_{max}$ is chosen
   In  order to improve the cancellation of unwanted term   for the difference in $\Delta t$ interval, improving so convergence, we could choose (see (\ref{Estremi}) ): 

  \be \label {0Codition}
  \Delta t \times \ln(\rho) = \frac {2 \pi}{\ln(p^*)}  \ln(\rho)   = \mu \times 2 \pi  \ ; \ with \  \mu=  m \in N \ , \ i.e. \  \rho=(p^*)^m
  \ee

 \noindent     In a simplistic view we could say:                                                                                                             
    
    \noindent $R(t,\epsilon,p_{max})$ is mainly an oscillation  with period close to $ \frac{2 \pi}{\ln(p_{max})}$, if we take $\Delta t= m \times \frac{2 \pi}{\ln(p_{max})} \ ; \ m \in N$  we filter away  the main part of $R$ for each  $  \rho=  p_{max} \times r $, when  $\frac{\ln(r)}{\ln(p_{max} ) } $ is small. So we have cancellation of main part of $R$} 
     not only for $\epsilon> 1/2$, but, also for  $\epsilon> 0$

  \begin{comment}  
     
       In general the expressions:
  
  $$
   \ln \left[    \prod_{p=2}^{p=p_{max}}         (1-p^{-s})          \right]
  -\int_2^{p_{max}}x^{-s}dLi(x)   \ \ ;  \forall t \ 
  $$ 
  
  or  (see (\ref {Estremi})  )
  
    $$ \frac  {\Delta \left \{
   \ln \left[    \prod_{p=2}^{p=p_{max}}         (1-p^{-s})          \right]
  -\int_2^{p_{max}}x^{-s}dLi(x) \right\}  } {\Delta t}  \ \ ;  \forall t \  and  \ \forall p^* : \ with \ \Delta t = \frac{2 \pi}{\ln(p^*)}$$

   \noindent  are absolutely  convergent, for $p_{max} =(p^*)^m \rightarrow \infty$ \footnote {This relation between $ p^*$ and $p_{max}$ will be clear in the following},  not only for $\epsilon> 1/2$, but, also for  $\epsilon> 0$

        \noindent  This does not mean (until now) that their value can be used for    $\ln [\zeta(s)(s-1)]$ or $\frac{\Delta \angle{\zeta(s)(s-1)}}{\Delta t} $  also for $\epsilon>0$, because  of condition  $\Re(s)>1$ in (\ref{1p1}). To attain this result we need next step.

  \end{comment}

\end{document}